\documentclass[11pt,a4paper]{article}
\usepackage{amssymb,amsmath}
\date{}
\author{}
\title{}

\title{Left invariant semi Riemannian metrics on quadratic Lie groups}
\author{A.Medina  \and  Ph. Revoy}

\date{}

\begin{document}

\title{\LARGE Lattices in some Symplectic or Affine Nilpotent Lie groups}

\maketitle

\vspace{.3in} \noindent Alberto Medina, Philippe Revoy\\D\'epartement des Math\'ematiques, C.C. 051\\%
Universit\'e Montpellier 2 , UMR CNRS 5149\\ Place E. Bataillon, 34095 
Montpellier cedex 5, France \\ medina@math.univ-montp2.fr \\revoyph@gmail.com



\abstract{ The main aim of this paper is the description of a large class of lattices in some nilpotent Lie groups, sometimes filiforms, carrying  a flat left invariant linear connection and often a left invariant symplectic form. As a consequence we obtain an infinity of, non homemorphic, compact affine or symplectic nilmanifolds. 
 We review some new facts about the geometry of compact symplectic nilmanifolds and we describe symplectic reduction for these manifolds. For the Heisenberg-Lie group, defined over a local associative and commutative finite dimensional real algebra, a necessary and sufficient condition for the existence of a left invariant symplectic form, is given. Finally in the symplectic case we show that a lattice in the group determines naturally lattices in the  double Lie group corresponding to any solution of the classical Yang-Baxter equation.  

\section{Introduction}

An affine manifold is a differential manifold $M$ with a special atlas of coordinate charts such that the coordinate changes extend to affine automorphisms of $ \mathbb{R}^{n}$. These distinguished charts are called affine charts. Given an affine structure on $M$ is equivalent to have a flat and torsion free linear connection $\nabla $ on $M$. A tensor field on $M$, for example a Poisson bivector, is called polynomial if in affine coordinates its coefficients are polynomial functions. For some results on affine geometry and polynomial tensor fields see ( [1],[2],[3],[4],[5] ). If $M=G$ is a Lie group and $\nabla $ is a left invariant affine structure on $G$ the pair $(G,\nabla)$ is called an\textbf{ affine Lie group}. From the infinitesimal point of view this means that the Lie algebra $ {\cal G}=T_{\epsilon}G$ of $G$ is endowed with a left symmetric product $ab:=L_{a}b=R_{b}a$ compatible with the Lie bracket of ${\cal G}$ in the sense that we have $ab-ba=[a,b]$. A pair $(G,\omega^{+})$ where $\omega^{+}$ is a left invariant symplectic form on $G$ is called a \textbf{symplectic Lie group}. For elements of the theory of symplectic Lie groups see [6]. The authors have described in [7],( see also [8],[9] ), the Lie algebras corresponding to nilpotent symplectic Lie groups in terms of the so-called symplectic double extension for the Lie algebras. For the structure of symplectic Lie groups see [10]. It is well known, see [11], that the formula 
\begin{center}
\begin{equation}
 \omega(ab,c)=-\omega(b,[a,c])
\end{equation}
\end {center}   
with $\omega^{+}_{\epsilon}=\omega$ and $\epsilon$ the unit in $G$, defines a left invariant affine structure on $G$. These two structures are used in [12] to prove the existence of polynomial Poisson tensors on some compact solvmanifolds. 

A pair $(G,\alpha^{+})$ where $\alpha^{+}$ is a left invariant contact form on the Lie group $G$, is called a \textbf{contact Lie group.} There are very strong links between the symplectic and the contact filiform Lie groups ([13]).

If a compact nilmanifold $M:=\Gamma\backslash G$ is endowed  with a symplectic form $\Omega$ , it is well known that there exists a left invariant symplectic form $\omega^{+}$ on $G$ inducing a symplectic form on $M$ cohomologous to $\Omega$ ([14]).The study of the geometry of  compact symplectic nilmanifolds $(M,\Omega)$ with $\Omega$ induced by $\omega^{+}$, called here \textbf{Symplectic Nilmanifolds}, is the other aim of this work. For a description of lattices in simply connected 4-dimensional solvable symplectic Lie groups, see [15].

\section{Lattices in 6-dimensional real 2-nilpotent symplectic Lie groups}
\subsection{ Case with derived Lie group 2-dimensional}
In this subsection ${G}$ is a 6-dimensional real 2-step nilpotent symplectic Lie group with Lie algebra ${\cal G}$ and derived group of dimension two. Denote by ${\cal H}_{1}(\mathbb{A})$ the 3-dimensional Heisenberg Lie $\mathbb{A}$-algebra where $\mathbb{A}$ is the field of real numbers, the field of complex numbers, or the ring $\mathbb{D}$ of dual numbers. The corresponding real Lie algebra ${\cal G}$ is one from ${\cal H}_{1}(\mathbb{C})$,  ${\cal H}_{1}(\mathbb{R}){\times}{\cal H}_{1}(\mathbb{R})$ or ${\cal H}_{1}(\mathbb{D})$, as shown below.  Each of these algebras is symplectic i.e. have an invertible scalar 2-cocycle (see [13]). In this section we want to study the lattices in $G$.


\subsubsection {Real Lie algebra structures on $\cal G $}

Let $A$ be a 6-dimensional 2-nilpotent real Lie algebra with center $ Z(A)=D(A)$ of dimension two. If $V$ is a vector subspace of $A$ such that $A=V \oplus D(A)$ , the bracket on $A$ is given by an onto linear map $ u: \Lambda^{2}(V)\rightarrow D(A)$. Consequently the transpose map of $ u $ is injective. Hence, given  $A$ is given  a vector subspace $W$ of dimension two of $\Lambda^{2}(V^{*})$ with maximal support i.e. $ \bigcap_{\eta\in W}Ker(\eta)  =0$ (*).

The Pfaffian on $\Lambda^{2}(V^{*})$ induces a quadratic map $ q:W\rightarrow \Lambda^{4}(V^{*})$, $\eta\mapsto Pf\eta:=\frac{1}{2}\eta^{2}$. If $ q $ were null there would be a basis $\{{\eta_{1},\eta_{2}}\}$ of $ W $ such that $\eta_{1}^{2}=\eta_{2}^{2}=\eta_{1}\wedge\eta_{2}= 0$. As $Ker\eta_{i}=P_{i}$ is a plane, $ \eta_{1}\wedge\eta_{2}= 0$ implies that $ P_{1}\bigcap P_{2}=D $ is a line and (*) is not verified. Now, any non zero element of $\Lambda^{4}(V^{*})$ determines, via $q$, a quadratic form on $ W $ we still note $q$. There are three cases for $q$  : 1) $q$ has rank two and is isotropic , 2) $q$ has rank two and is anisotropic, 3) $q$ has rank one.
 
In the first case there exists a basis $\{{\eta_{1},\eta_{2}}\}$ of $W$ such that $\eta_{1}^{2}=0=\eta_{2}^{2}$ and $ \eta_{1}\wedge\eta_{2} \neq 0$. If $\eta_{1}=f_{1}\wedge g_{1} $ and $\eta_{2}=f_{2}\wedge g_{2} $ then $\{{f_{1},g_{1},f_{2},g_{2}\}}$ is a basis of $V^{*}$. If  $\{{e_{1},e_{2},e_{3},e_{4}\}}$ denote its dual basis ,we have $\eta_{1}=e_{1}^{*}\wedge e_{2}^{*}$ and $\eta_{2}=e_{3}^{*}\wedge e_{4}^{*}$. Hence, the Lie algebra $A=V \oplus D(A)$ has two non null brackets $[e_{1},e_{2}]=:e_{5}$ and $[e_{3},e_{4}]=:e_{6}$. Moreover as $D(A)$ is a plane, the Lie algebra $A$ is isomorphic to the direct product ${\cal H}_{1}(\mathbb{R}){\times}{\cal H}_{1}(\mathbb{R})$.

Suppose $q$ anisotropic with rank two. Let $A_{\mathbb{C}}$ , $W_{\mathbb{C}}$ and  $q_{\mathbb{C}}$ be the corresponding complexifications. It is clear that $q_{\mathbb{C}}$ is isotropic and there exists $\eta\in W_{\mathbb{C}}$ non zero such that $\eta^{2}=0$. Putting $\eta=\eta_{1}+i\eta_{2}$ with $\eta_{i}\in W $ , we have $\eta_{1}^{2}=\eta_{2}^{2}\neq 0 $ and $ \eta_{1}\wedge\eta_{2}=0$. As $\eta^{2}=0$ then $\eta=f\wedge g $ with $f:=f_{1}+i f_{2}$ and $g=g_{1}+i g_{2}$ where $f_{1},f_{2},g_{1},g_{2}\in V^{*}$. The fact that $ \eta_{1}=Re\eta,  \eta_{2}=Im\eta $ implies that
$ \eta_{1}= f_{1}\wedge g_{1}-f_{2}\wedge g_{2}$ ,  $ \eta_{2}= f_{1}\wedge g_{2}+f_{2}\wedge g_{1}$ and $B^{*} =\{{f_{1},g_{2},f_{2},g_{1}\}}$ is a basis of $V^{*}$. Let $\{{e_{1},e_{2},e_{3},e_{4}\}}$ be the dual basis of $B^{*}$. We have then, $ \eta_{2}=e_{1}^{*}\wedge e_{2}^{*} + e_{3}^{*}\wedge e_{4}^{*}$ ; $ \eta_{1}=e_{1}^{*}\wedge e_{4}^{*} + e_{2}^{*}\wedge e_{3}^{*}$ . Consequently the product in the Lie algebra $A=V \oplus D(A)$ is given by the brackets $ [e_{1}, e_{2}]=[e_{3}, e_{4}]=:e_{5}$ and $ [e_{1}, e_{4}]=[e_{2}, e_{3}]=:e_{6} $. Now consider a basis $\{{f_{1},f_{2},f_{3}\}}$ of the Heisenberg complex Lie algebra of dimension three with $[f_{1},f_{2}]=f_{3}$. The real basis $\{{f_{1},if_{2},if_{1},f_{2},if_{3},f_{3}\}}$ denoted by $\{{e_{1},e_{2},e_{3},e_{4},e_{5},e_{6}\}}$ verifies $ [e_{1}, e_{2}]=[e_{3}, e_{4}]=:e_{5}$ and $ [e_{1}, e_{4}]=[e_{2}, e_{3}]=:e_{6} $ , so $A$ is isomorphic to ${\cal H}_{1}(\mathbb{C})$. 

Finally, if $ q $ has rank 1, there is a basis $ \{\eta_{1},\eta_{2}\} $ of $ W $ such that $\eta_{1}^{2}\neq o $ ; $\eta_{2}^{2}=\eta_{1}\wedge\eta_{2}= 0 $. From $\eta_{2}^{2}= o $, it comes, $\eta_{2}=f\wedge g $ and $ \eta_{1}\wedge\eta_{2}= 0 $ implies $\eta_{1}= f \wedge h + g\wedge k $ with  $B^ {*} =\{f,g,h,k\}$  a basis of $ V^{*} $.  Let $\{{e_{1},e_{2},e_{3},e_{4}}\} $ be the dual basis of $B^{*}$. We have $ \eta_{1}=e_{1}^{*}\wedge e_{3}^{*} + e_{2}^{*}\wedge e_{4}^{*}$ and $\eta_{2}=e_{1}^{*}\wedge e_{2}^{*}$. Then, in terms of $ B $ the bracket in $ A $ is given by $ [e_{1},e_{3}]= [e_{2},e_{4}]=:e_{5}$ and $[e_{1},e_{2}]= :e_{6} $. To show that the Lie algebra $A$ is isomorphic to ${\cal H}(\mathbb{D})$, consider a $\mathbb{R}[\epsilon]$-basis $\{{f_{1},f_{2},f_{3}}\} $ of the last algebra with $[f_{1},f_{2}]=f_{3}$. As $ \epsilon^{2}=0$ in the $\mathbb{R}$-basis $\{e_{i}; 1\leq i\leq 6\}=\{f_{1},f_{2},\epsilon f_{2},-\epsilon f_{1},\epsilon f_{3},f_{3}\} $ we have $[e_{1},e_{2}]=e_{6}$  ; $[e_{1},e_{3}]= e_{5}=[e_{2},e_{4}]$  ; $[e_{3},e_{4}]=0$ and the two real Lie algebras $A$ and ${\cal H}_{1}(\mathbb{D})$ are isomorphic.

\subsubsection {Rational Lie algebra structures on ${\cal G}$ and lattices}

Let $A$ be a rational Lie algebra of dimension 6 such that the real Lie algebra  $ A\otimes \mathbb{R} $ obtained by extension of the base field to $\mathbb{R}$ is isomorphic to one of Lie algebras described in 2.1. It is clear that $A$ is 2-step nilpotent and $ Z(A)=D(A)$ has dimension two. The choice of a subspace $V$ of $A$ as above determines a plane $W$ of $\Lambda^{2}(V^{*})$ and a quadratic map  $ q:W\rightarrow \Lambda^{4}(V^{*})$. 
If the rank of $q$ is 1, there is a basis $ \{\eta_{1},\eta_{2}\} $ of $ W $ such that 
$\eta_{1}^{2}\neq o $ ; $\eta_{2}^{2}=\eta_{1}\wedge\eta_{2}= 0 $ and the above calculation  shows the existence of a basis of $A$ in which the multiplication table is that of ${\cal H}_{1}(\mathbb{D})$. Hence, the corresponding real 1-connected Lie group $G$ is the unimodular group of matrix 

\begin{center}

$\begin{pmatrix} 
 
1 & u & v \\
0 & 1 & w \\
0 & 0 & 1

\end{pmatrix}$

\end{center}
\begin{flushleft}
with $u,v,w \in \mathbb{R}[\epsilon]=:\mathbb{D}$.\\
\end{flushleft}
 In this case, any lattice in $G$ is commensurable to the lattice $\Gamma$ whose elements are the matrices in $G$ with $u,v,w\in \mathbb{Z}[\epsilon]$.

If $q$ has rank two, the choice of a basis of $V$ determines a quadratic form on $W$. Moreover, it is cleat that the quadratic form produced for another basis choice, is collinear with the first. In other words, the set of the $q$ maps obtained is in bijection with the $\mathbb{Q}^{*}$- orbits in the set of quadratic forms of rank two on $W$ where the action is given by homotheties.
In every equivalent class there is an unique representative of the form $X^{2}+dY^{2}$ with $d$ the common discriminant of the quadratic forms of the class. It is clear that $d$ can be chosen a squarefree integer.  

Hence we have :

If $d< 0$, the real algebra $ A\otimes \mathbb{R} $ is isomorphic to ${\cal H}_{1}(\mathbb{R}){\times}{\cal H}_{1}(\mathbb{R})$ .

If $d > 0$, $ A\otimes \mathbb{R} $  is isomorphic to  ${\cal H}_{1}(\mathbb{C})$ .

If $ d\neq -1 $ then ${\cal K }:= \mathbb{Q}(\sqrt {-d} )$ is a field, where $ \mathbb{Q} $ is the field of rational numbers. As a consequence, the extension of the scalar field to ${ \cal K}$ produces the Lie algebra direct product ${{\cal H}}_{1}({\cal K})\times {\cal H }_{1}({\cal K})$. Consider  $q_{A}:W\rightarrow \Lambda^{4}(V^{*})$, the form composed of the injection of $W$ in $\Lambda^{2}(V^{*})$ with $\gamma_{2}=Pf$ where $Pf$ is the Pfaffian. Obviously $q_{A}$ is proportional to $X^{2}+dY^{2}$ and its extension to the tensor product $A\otimes \mathbb{Q}(\sqrt {-d})$, relative to $\mathbb{Q}$, becomes isotropic: that is, there exists $ u\in W\otimes{\cal K }$, $u\neq0 $,  $\gamma_{2}(u)=0$ and $u$ is decomposable i.e. $u= (f_{1}+\sqrt {-d}f_{2})\wedge (f_{3}+\sqrt {-d}f_{4})$ with $f_{j}\in V^{*}$. As $u=( f_{1}\wedge f_{3}-d f_{2}\wedge f_{4})+\sqrt {-d}( f_{1}\wedge f_{4}+ f_{2}\wedge f_{3})=u_{1}+\sqrt {-d}u_{2}$ where $u_{1}=\frac{u+\overline{u}}{2}$, $u_{2}=\frac{u-\overline{u}}{2 \sqrt {-d}}$ with $\overline{u}$ the conjugate of $u$ .Since $u_{1}$ and $u_{2}$ have rank 4 , $B^{*}=\{f_{1},f_{2},f_{3},f_{4}\}$ is a basis of $V^{*}$ and $\{u_{1},u_{2}\}$ is a basis of $W$ over $\mathbb{Q}$ . At this point we can write the bracket of the Lie algebra $A$ over $\mathbb{Q}$ .Let $ B= \{e_{1},e_{2},e_{3},e_{4}\}$ the dual basis of $B^{*}$. As $ u_{1}=e_{1}^{*}\wedge e_{3}^{*}-d e_{2}^{*}\wedge e_{4}^{*} $ and $u_{2}=e_{1}^{*}\wedge e_{4}^{*}+ e_{2}^{*}\wedge e_{3}^{*} $, the bracket in $A$ is given by 

\begin{center}
\begin{equation}
[e_{1},e_{4}]=[e_{2},e_{3}]=e_{5} , [e_{1},e_{3}]=e_{6} ,  [e_{2},e_{4}]=-d e_{6}
\label{eq:ED1}
\end{equation}
\end{center}

\begin{flushleft}

where $d$ is a non null squarefree integer \\
\end{flushleft}

The formula (2) gives all the $\mathbb{Q}$ -structures of Lie algebra on the real Lie algebras 
${\cal H}_{1}(\mathbb{R}){\times}{\cal H}_{1}(\mathbb{R})$ and ${\cal H}_{1}(\mathbb{C})$ if $d>0$ or $d<0$ respectively (see [16], page 47).

We want, in the following,  describe a lattice in each class of commensurability.

Notice if $d > 0$ , the corresponding 1-connected Lie group is the 3-dimensional Heisenberg Lie group  $G= H_{1}(\mathbb{C})$ and we have the following tower of groups $G= H(\mathbb{C}) \supset G_{d}= H({\cal K})\supset\Gamma_{d}=H({\cal O}_{-d})$   where ${\cal O}_{-d}$ is the ring of integers  of ${\cal K }$. In other words:
\begin{center}
 
${\cal O}_{-d} =\mathbb{Z}(\sqrt {-d})$ if $ d \equiv 1,2 ( mod\ 4)$  

${\cal O}_{-d}=\mathbb{Z}[\frac{1+\sqrt{-d}}{2}]$ if  $d\equiv 3 ( mod \  4 )$

\end{center}

\textbf{Remark 1.} Obviously the lattice $\Gamma_{d}:=N_{d}(\mathbb{Z})$ is embedded into $N_{d}(\mathbb{Q})=: H_{1}({\cal K})$ . Also if $d_{1}\neq d_{2}$ are two squarefree positive integers, the lattices $N_{d_{1}}(\mathbb{Z})$ and $N_{d_{2}}(\mathbb{Z})$ are not commensurables.

From the multiplication in the group $G'$ of matrices 
\begin{center}

$\begin{pmatrix}

1 & a & c \\
0 & 1 & b \\
0 & 0 & 1

\end{pmatrix}$

\end{center} 
\begin{flushleft}
with $a,b,c \in {\cal O}_{-d}$ , we can deduce the product in the lattice $H_{1}({\cal O}_{-d})$. In fact, we have $N_{d}(\mathbb{Q})=\mathbb{Q}^{6}$ and the product in the lattice $N_{d}(\mathbb{Z})$ is induced by the product of $G'$. More precisely, if $ d\equiv 1,2 (mod \  4)$ putting $a=x_{1}+\varpi x_{2}$ and   $b'=y_{3}+\varpi y_{4}$ with $\varpi^{2}=-d $ it is found that the the products in the groups $H_{1}(\mathbb{Z}[\sqrt{-d}])$ and $N_{d}(\mathbb{Z})$ are identical.
\end{flushleft}

In a similar way, the equality $\varpi^{2}=\varpi-\frac{1+d}{4}$ determines the product in the lattice when $d\equiv 3 (mod\  4 )$.

 Let $d\in\mathbb{Z}$ be squarefree . If $d > 0 $ then $\mathbb{Q}( \sqrt{-d})$ is a subfield of $ \mathbb{C}$ whereas if $d<0$ , $\mathbb{Q}( \sqrt{-d})$ is a subfield of $\mathbb{R}\times \mathbb{R}$. In fact the map $\mathbb{Q}[\sqrt{-d}]=\frac{\mathbb{Q}[X]}{(X^{2}+d)}\rightarrow \mathbb{R}\times \mathbb{R}$ , $\overline{\alpha+\beta X }\mapsto (\alpha+\sqrt{-d} \beta,\alpha-\sqrt{-d} \beta)$ is an injective homomorphism of rings. Hence $H_{1}({\cal O}_{-d})\subset H_{1}({\cal K})$ and the last one is included in $H_{1}(\mathbb{C})$ or in $H_{1}(\mathbb{R}\times \mathbb{R})$ according to the sign of $d$.

Recall that a field of the form $\frac{\mathbb{Q}[X]}{(X^{2}+d)}$ where $d$ is an integer and $d\neq 0,-1 $ is called a quadratic field. If $d < 0$ the field is said real, if $d>0$, the field is imaginary.

In the these terms we have shown the following result:

\textbf{Theorem 1.} {a. Each  nilpotent real Lie algebras ${\cal H}_{1}(\mathbb{C})$ and ${\cal H}_{1}(\mathbb{R}){\times}{\cal H}_{1}(\mathbb{R})$ has an infinity of $\mathbb{Q}$ -Lie algebra structures  corresponding to squarefree integers $d$. Consequently the corresponding simply connected real Lie groups $ H_{1}(\mathbb{C})$ and $H_{1}(\mathbb{R}){\times}H_{1}(\mathbb{R})$ have an infinity of commensurability class of lattices.

b. In $ H_{1}(\mathbb{C})$ each commensurability class contains a lattice isomorphic to $H_{1}({\cal A})$ where ${\cal A}$ is the ring of integers of the imaginary quadratic field $\mathbb{Q}(\sqrt{-d})$. In $H_{1}(\mathbb{R}){\times}H_{1}(\mathbb{R})$ every class of commensurability either contains a lattice isomorphic to $H_{1}({\cal A})$ with ${\cal A}$ the ring of integers of the real quadratic field $\mathbb{Q}(\sqrt{-d})$ if $d < -1$ or $H_{1}(\mathbb Z)\times H_{1}(\mathbb Z)$ if $d=-1$.}

\textbf{Remark 2.} {In terms of a basis $\{e_{i}, 1\leq i\leq 6 \}$ the brackets in the $\mathbb Q $-Lie algebra, as $ d>1 $ below, are given by 
\begin{center}
\begin{equation}
[e_{1},e_{2}]=[e_{3},e_{4}]=e_{5} , [e_{1},e_{3}]=e_{6} ,  [e_{2},e_{4}]= -d e_{6}
\end{equation}
\end{center}
where $d$ is a squarefree integer.
Moreover, the product $\sigma\times \tau $  in the  algebraic group $N_{d}(\mathbb Z )$ for $\sigma=(x_{i})$ and $\tau=(y_{i})$ is equal to either
\begin{equation}
\sigma\tau=(x_{1}+y_{1},..,x_{4}+y_{4}, x_{5}+y_{5}+x_{1}y_{3}+dx_{2}y_{4},x_{6}+y_{6}+x_{1}y_{4}+x_{2}y_{3})
\label{eq:ED1}
\end{equation}  
or 
\begin{equation}
\sigma\tau=(x_{1}+y_{1},..,x_{4}+y_{4}, x_{5}+y_{5}+x_{1}y_{3}+(\frac{d-1}{4})x_{2}y_{4},x_{6}+y_{6}+x_{1}y_{4}+x_{2}y_{3}+x_{2}y_{4})
\end{equation}

\begin{flushleft}
if $d\equiv 2,3  ( mod\  4 )$ or $d\equiv 1  ( mod\  4 )$ respectively ).
\end{flushleft}

\subsection {Lattices in the cotangent bundle of the real Heisenberg-Lie group}  

Denote by $G:= T^{*} H_{1}(\mathbb{R}) $ the simply connected Lie group  whose Lie algebra ${\cal G}$ is the semi-direct product of the real Heisenberg Lie algebra ${\cal H}_{1}$ by the abelian Lie algebra ${\cal H}_{1}^{*}$ via the coadjoint action of ${\cal H}_{1}$. We will describe the lattices of $G$. This group is a quadratic  2-step nilpotent Lie group with bi-invariant metric of type $ (3,3)$. Also it carries a left invariant symplectic form: the corresponding scalar 2-cocycle $\Omega$ is defined by an invertible derivation $\delta=(\lambda,-\lambda^{t}) $ where $\lambda$ is an invertible derivation of ${\cal H}$.

It is clear that the Lie algebra ${\cal G}$ is isomorphic to the vector space $\mathbb{R}^{3}\oplus \Lambda^{2}(\mathbb{R}^{3})$ endowed with the Lie bracket $[(x,u),(y,v)]=(0,x \wedge y)$ and that ${\cal G}$ has a unique rational Lie algebra structure given by this bracket. So, any two lattices in $G$ are commensurable.

It is useful, for our purpose, to identify $G$ and the vector space $\mathbb{R}^{6}$ endowed with the product given for $\sigma=(x_{1},x_{2},x_{3},y_{1},y_{2},y_{3})$ and $\tau= (x'_{1},x'_{2},x'_{3},y'_{1},y'_{2},y'_{3})$  as 
\begin{equation}
\sigma\tau=(x_{1}+x'_{1}+y_{2}y'_{3},x_{2}+x'_{2}+y_{3}y'_{1},x_{3}+x'_{3}+y_{1}y'_{2},y_{j}+y'_{j})
\end{equation}
where $1\leq j\leq 3$.

The subset $\Gamma_{0}$ of $G$ whose  elements have integer components is a lattice of $G$ and any lattice of $G$ is commensurable to $\Gamma_{0}$. Moreover its derived group is equal to its center. We have :
 
\textbf{Theorem 2.} Up to isomorphism a lattice $\Gamma $ of $ G $, as above , is characterized  by three positive integers $d_{1}, d_{2},d_{3}$ such that $1\leq d_{1}|d_{2}|d_{3}$. More precisely, $\Gamma $ is isomorphic to the group 
\begin{center}
$\Gamma_{d_{1},d_{2},d_{3}}=\{(a_{i},b_{j});a_{i},b_{j}\in\mathbb{Z},1\leq i,j\leq3\}$
\end{center}
with multiplication given by
\begin{center}
$(a_{i},b_{j})(a'_{i},b'_{j})=(a_{1}+a'_{1}+d_{1}b_{2}b'_{3},a_{2}+a'_{2}+d_{2}b_{3}b'_{1},a_{3}+a'_{3}+d_{3}b_{1}b'_{2},b_{j}+b'_{j})$
\end{center}   

\textbf{Proof}. A calculation shows that $\Gamma_{d_{1},d_{2},d_{3}}$ is a group. It is a subgroup of $G$ because the extension of its multiplication to $\mathbb{R}^{6}$ gives a group isomorphic to $G$. Also it is clear that $\Gamma_{d_{1},d_{2},d_{3}}$ is a lattice of $G$ , so commensurable to  $\Gamma_{0}$.

In the other hand, the center and the derived group of $\Gamma_{d_{1},d_{2},d_{3}}$ are respectively 
\begin{center}
$ Z(\Gamma_{d_{1},d_{2},d_{3}})=\{( (a_{i}),0,0,0);a_{i}\in\mathbb{Z}\}$ 

$D(\Gamma_{d_{1},d_{2},d_{3}})=\{(d_{1}a_{1},d_{2}a_{2},d_{3}a_{3},0,0,0);a_{i}\in\mathbb{Z}\}$

\end{center} 
\begin{flushleft}
Hence, the factor group  $Z/D$  is isomorphic to the direct product of  three cyclic groups $C_{d_{1}}\times C_{d_{2}}\times C_{d_{3}}$. Consequently different triples of integers, as in the Theorem 2, produce non isomorphic lattices.
\end{flushleft}
The following presentation of the group $\Gamma_{d_{1},d_{2},d_{3}}$ by  generators and relations, will be useful for the proof. Denote by $z_{1},z_{2},z_{3},y_{1},y_{2},y_{3}$ the elements of $\Gamma_{d_{1},d_{2},d_{3}}$ with only one non zero component which is equal to $1$ in the place $i$ for $z_{i}$ and $j+3$ for $y_{j}$. The $z_{i}$ generate the center  and satisfy the relations:
\begin{center}
$z_{k}z_{l}=z_{l}z_{k}$ and $ y_{j}z_{i}=z_{i}y_{j}$ for any $i,j,k,l\in\{1,2,3\}$
\end{center}
and a direct calculation gives the other relations
\begin{center}
\begin{equation}
y_{2}y_{3}=y_{3}y_{2}z_{1}^{d_{1}} , y_{3}y_{1}=y_{1}y_{3}z_{2}^{d_{2}} ,y_{1}y_{2}=y_{2}y_{1}z_{3}^{d_{3}}.
\end{equation}
\end{center}
  
Conversely, a group generated by six elements satisfying these relations with the $d_{i}$ as above is isomorphic to $\Gamma_{d_{1},d_{2},d_{3}}$.

Finally we show that any lattice $\Gamma $ of $G$ has such a presentation.

It is clear that a lattice $\Gamma $ of $G$ is a 2-step nilpotent group: its center $Z$ contains its derived group $D$. Also $D$ is a sublattice of the lattice $Z$  of the 3-dimensional real vector space $Z(G)$. Hence $Z/D$ is a group direct product of three cyclic groups $C_{d_{1}},C_{d_{2}}$ and $ C_{d_{3}}$ where the $d_{i}$ are positive integers and 
 $1\leq d_{1}|d_{2}|d_{3}$ . Then, using  multiplicative notation , there is a basis $\{z_{1},z_{2},z_{3}\}$ of the free Abelian group $Z$ such that $\{z_{1}^{{d}_{1}},z_{2}^{{d}_{2}},z_{3}^{{d}_{3}}\}$ is a basis of $D$.
 
Let $K$ be an arbitrary group and $ \left\langle g,h\right\rangle:=g^{-1}h^{-1}gh$ the commutator of $g$ and $ h\in K $ . It is easy to verify the formula 
\begin{center}
$ \left\langle g,h_{1}h_{2}\right\rangle= \left\langle g,h_{2}\right\rangle \left\langle g,h_{1}\right\rangle^{h_{2}}$ 
\end{center}
where $a^{b}:=b^{-1}ab$.

Applying this formula to $\Gamma $ we get :
\begin{center}
$\left\langle x,yz\right\rangle=\left\langle x,y\right\rangle\left\langle x,z\right\rangle$ 
\end{center}
because two commutators commute and $D$ is in the center of $\Gamma $. 

Now the commutator map 
\begin{center}
$\Gamma\times \Gamma\rightarrow D$ , $(x,y)\rightarrow \left\langle x,y\right\rangle $
\end{center}
induces an antisymmetric  map $\varphi: \Gamma/Z \times \Gamma/Z\rightarrow D $. But the last formula shows that $\varphi $ is additive in the second variable. Consequently it is bi-additive and then determines an homomorphism of Abelian groups $\overline{\varphi}:\Lambda^{2}( \Gamma/Z  )\rightarrow D $. This homomorphism is onto. Moreover , $ \Gamma/Z $ can be identifyed with a lattice of the real space $G/Z(G)$. As $D$, $\Gamma/Z$ and $\Lambda^{2}( \Gamma/Z )$ are free Abelian groups of finite type and rank 3, it follows that $\overline{\varphi}$ is an isomorphism of groups. Denote by $\{u_{1},u_{2},u_{3}\}$
the basis of $\Lambda^{2}( \Gamma/Z )$ such that $\overline{\varphi}(u_{i})=z_{i}^{{d}_{i}}$. As $\mathbb{Z}$ is a principal ring, there exists a basis $\{t_{1},t_{2},t_{3}\}$ of $\Gamma/Z$ such that  $u_{1}=t_{2}\wedge t_{3} $, $u_{2}=t_{3}\wedge t_{1} $ and $u_{3}=t_{1}\wedge t_{2} $ . Taking $y_{i}\in\Gamma $ with  $\overline {y_{i}}=t_{i}$ for $1\leq i\leq 3$ we have 
$ \left\langle y_{2},y_{3}\right\rangle=z_{1}^{{d}_{1}}$....Hence it is clear that $\Gamma $ is generated by $\{y_{1},y_{2},y_{3},z_{1},z_{2},z_{3}\}$ with the relations:
\begin{center}
$y_{2}^{-1}y_{3}^{-1}y_{2}y_{3}=z_{1}^{{d}_{1}}$,

$y_{3}^{-1}y_{1}^{-1}y_{3}y_{1}=z_{2}^{{d}_{2}}$,

$y_{1}^{-1}y_{2}^{-1}y_{2}y_{1}=z_{3}^{{d}_{3}}$.

$z_{i}x=xz_{i}$ 
\end{center}
\begin{flushleft}
for every $x \in \Gamma $.
\end{flushleft}
So, $\Gamma $ is isomorphic to $\Gamma_{d_{1},d_{2},d_{3}}$.

 \hfill $\square$\\

\section{Lattices in some filiform Lie groups}

For an integer $n\geq 2$ let ${\cal L}_{n}(\mathbb{K})=Span_{\mathbb{K}}\{e_{0},e_{1},...,e_{n}\}$ be the  filiform Lie algebra over the field ${\mathbb{K}}$ with the only nonzero brackets given by $ [e_{0},e_{i}]=e_{i+1}$ for $i=1,2,...,n-1$. As ${\cal L}_{n}(\mathbb{K})$ has an invertible derivation,its bracket is compatible with a left symmetric product.

\subsection {Rational Lie algebra structures in ${\cal L}_{n}(\mathbb{R})$} 

The  algebras ${\cal L}_{n}(\mathbb{K})$ are the only filiform Lie algebras having a unique Abelian ideal of codimension $1$. In fact a codimension $1$ ideal $I$ contains the derived ideal $J=Span_{\mathbb{K}}\{e_{2},...,e_{n}\}$ and is determined by a line in the quotient algebra by $J$, so $I=Span_{\mathbb{K}}\{\alpha e_{0}+\beta e_{1},e_{2},...,e_{n}\}$ with $(\alpha,\beta)\neq (0,0)$. If $\alpha\neq 0$, $I$ is not Abelian whereas, if $\alpha=0$, $I=Span_{\mathbb{K}}\{e_{1},...,e_{n}\}$ is the unique Abelian codimension $1$ ideal of ${\cal L}_{n}(\mathbb{K})$. We have :

\textbf{Lemma 1.} Let $k\subset\mathbb{K}$ be a subfield of the field $\mathbb{K}$ and ${\cal G}$ a $k$-Lie algebra such that the Lie algebra ${\cal G}\otimes_{k}\mathbb{K}$ is isomorphic to ${\cal L}_{n}(\mathbb{K})$. Then ${\cal G}$ has a unique Abelian ideal of codimension 1 and  is isomorphic to ${\cal L}_{n}(k)$.

\textbf{Proof.} It is clear that $D({\cal G}){\otimes}_{k}\mathbb{K}$ is isomorphic to 
$D({\cal L}_{n}(\mathbb{K}))=Span _{\mathbb{K}}\{e_{2},...e_{n}\}$ and hence $ D({\cal G})$ is a $k$-subspace of dimension $n-1$ of $ D({\cal L}_{n}(\mathbb{K}))$.

Consider a $k$-basis $\{f_{2},f_{3},...,f_{n}\}$ of $D({\cal G})$ completed by $\{f_{0},f_{1}\}$ in a $k$-basis $B$ of ${\cal G} $. As $B$ is also a $\mathbb{K}$-basis of ${\cal L}_{n}(\mathbb{K})$ we may write $e_{1}=\alpha f_{0}+\beta f_{1}+\sum_{i\geq 2} \gamma_{i}f_{i}$ with $\alpha,\beta,\gamma_{i}\in\mathbb{K}$. 

Also, as $ ad e_{1}|D({\cal L}_{n}(\mathbb{K}))=0$ and $ ad f_{i}|D({\cal L}_{n}(\mathbb{K}))= 0$ for $i\geq 2$ ,
\begin{center} 
$(\alpha ad f_{0}+\beta ad f_{1})|D({\cal L}_{n}(\mathbb{K}))=0$.
\end{center}
Now, if $ad f_{0}|D({\cal L}_{n}(\mathbb{K}))=0$ or $ad f_{1}|D({\cal L}_{n}(\mathbb{K}))=0$ then $D({\cal G})\oplus k f_{i}$ is an Abelian ideal of ${\cal G}$ of codimension 1.

Finally, suppose both adjoints are non zero. Using $ \{f_{2},...,f_{n}\} $ as a $\mathbb{K}$-basis of $D({\cal L}_{n}(\mathbb{K}))$, we get two $n-1$ square matrices, with coefficients in $k$, representing the restrictions of $ad f_{0}$ and $ad f_{1}$ to $D({\cal L}_{n}(\mathbb{K})$, which are  $\mathbb{K}$-linearly dependant. This implies that $\alpha$ and $\beta$ are $k$ linearly dependent i.e. $\alpha=\lambda\beta$ with $\lambda\in k-\{0\}$. Then we have $e_{1}=\beta (\lambda f_{0}+f_{1})+\sum_{i\geq 2}\gamma_{i}f_{i}$ with $\lambda f_{0}+f_{1}\in{\cal G}$. It is then clear that $Span_{k}\{\lambda f_{0}+f_{1},f_{2},...,f_{n}\}$ is a codimension $1$ ideal of 
${\cal G}$ and that ${\cal G}$ is isomorphic to ${\cal L}_{n}(k)$.  

 \hfill $\square$\\


Denote by $F_{n}$ the simply connected real Lie group with Lie algebra ${\cal L}_{n}(\mathbb{R})$. We have the following consequence of the Lemma 1 (see [16]).

\textbf{Lemma 2.} The Lie group $F_{n}$ has a unique commensurability class of lattices.

\subsection {Isomorphism classes of lattices in $F_{n}$ }

In the following, we want to describe the isomorphism classes of lattices in $F_{n}$. It is clear that our Lie group is isomorphic to the manifold ${\mathbb{R}}^{n}\times {\mathbb{R}}$ endowed with the product 
\begin{center}
$(v,t)(v',t')=(v+ (Exptd )(v'),t+t')$
\end{center}
where $d$ is a principal nilpotent endomorphism of the vector space ${\mathbb{R}}^{n}$. The center of the group is then $Z F_{n}=\{(v,0) ; v\in Ker d \}$ and we have the following canonical exact sequence of Lie groups
\begin{center}

$ 1\leftharpoonup Z F_{n}\rightarrow\stackrel {p_{n}}{F_{n}\rightarrow  F_{n-1}} \rightarrow 1 $

\end{center}

Obviously, the ascending  central sequence of $F_{n}$ is obtained by pullback by 
$p_{n}$ of the corresponding sequence of $F_{n-1}$. By induction on $n\geq2 $ , it follows that

\begin{center}
$0\subset C_{1}F_{n}\subset C_{2}F_{n}\subset...\subset C_{r}F_{n}\subset...$
\end{center}
is a sequence of Abelian groups and $ C_{r}F_{n}$ is isomorphic to ${\mathbb{R}}^{r}$ for $r\leq n-1$ and $C_{n}F_{n}=F_{n}$. Moreover as $C^{1}({\cal L}_{n})=C_{n-1}({\cal L}_{n})$ ; 
$C^{2}({\cal L}_{n})=C_{n-2}({\cal L}_{n})$ ...$C^{r}({\cal L}_{n})=C_{n-r}({\cal L}_{n})$ , we also know the descending central sequence of $F_{n}$.

It is clear that the set $\Gamma_{0}:={\mathbb{Z}^{n}}\times \mathbb{Z}$ endowed with the (semidirect) product given by

\begin{center}

$(a_{1},...,a_{n},k)(a'_{1},...,a'_{n},k')=((a_{1},...a_{n})+ g_{0}^{k}(a'_{1},...a'_{n}),k+k')$
\end {center}
with
\begin{center}
 $g_{0}^{k}(a'_{1},...a'_{n})= (a'_{1},a'_{2}+ka'_{1},a'_{3}+ka'_{2}+\binom{k}{2}a'_{1},...)$
\end {center}
\begin{flushleft}
is a lattice in $F_{n}$  and any lattice in $F_{n}$ is commensurable with $\Gamma_{0}$.\\
\end{flushleft}

Notice that the action of $\mathbb{Z}$ over \textsl{${\mathbb{Z}^{n}}$} is essentially given by the product of the matrix $g_{0}:=I_{n}+E_{2,1}+E_{3,2}+...E_{n,n-1}$ by the column vector $(a'_{1},a'_{2},...,a'_{n})^{T}$. Here $E_{i,j}$ is the n-square matrix with only one non zero coefficient , which is 1 in the i-line and the j-column.

In the following, we will describe the other lattices in $F_{n}$.

Let $\Gamma_{1}$ be a finite index subgroup of $\Gamma_{0}$ and consider the corresponding canonical sequence of groups: 
\begin {center}
$\Gamma_{1}\leftharpoonup \Gamma_{0}\rightarrow \mathbb{Z}$
\end{center} 
with  $\pi:\Gamma_{0}\rightarrow \mathbb{Z}$ ,$\pi(a_{1},...a_{n},k)=k$. It is clear that 
$\pi( \Gamma_{1})=d \mathbb{Z}$ where $d$ a strictly positive integer.

Let $\pi'$ be the restriction of $\pi $ to $\Gamma_{1}$. Obviously $L_{1}:=Ker \pi'\subset Ker\pi= {\mathbb{Z}^{n}}=\{(a_{1},...a_{n},0 );a_{i}\in \mathbb{Z} \}$ is a finite index subgroup of ${\mathbb{Z}^{n}}$ and we obtain the following canonical sequence of groups

\begin{center}
$ 1\leftharpoonup L_{1} \leftharpoonup \Gamma_{1}\rightarrow d  \mathbb{Z}\rightarrow 1 $

\end{center}

Thus  $\Gamma_{1}$  is the group semi-direct product of $L_{1}$ by $\mathbb{Z}$ with the action of $\mathbb{Z}$ over $L_{1}$ given by $g_{1}:=g_{0}^{d}|L_{1}$.

\textbf{ Let $\Gamma $ be a lattice of $F_{n}$} . This means that $\Gamma $ contains a sublattice of $\Gamma_{0}$ and this one is a finite index subgroup in $\Gamma$ and in $\Gamma_{0}$.

Let $H_{n}$ be the maximal Abelian normal subgroup ( of codimension 1 ) of $F_{n}$ : it is clear that ${\mathbb{Z}^{n}}$ and $L_{1}$ are lattices of $H_{n}$. Consider the canonical injection $j:\Gamma\leftharpoonup F_{n}$ and the canonical surjection $p:F_{n}\rightarrow F_{n}/H_{n}$ . Then we have the exact sequence of groups:

\begin{center}

$ 1\leftharpoonup H_{n}\leftharpoonup F_{n}\rightarrow  F_{n}/H_{n}\cong \mathbb{R}$

\end{center}

Since $(p\circ j)(\Gamma_{1})$ is a finite index subgroup of $(p\circ j)(\Gamma)$ , the last one is a lattice of $\mathbb{R}$. 

\textbf{Remark 3.} Consider the Lie group $G:=\mathbb{R}^{3}\oplus \Lambda^{2}\mathbb{R}^{3}$ , its lattice $\Gamma=\mathbb{Z}^{3}\oplus \Lambda^{2}\mathbb{Z}^{3}$ and  $H:=P\oplus \mathbb{R}^{3}$ where $P$ is any plane of $\mathbb{R}^{3}$. The factor group $G/H$ is isomorphic to $\mathbb{R}$ but  the image of $\Gamma$ by the canonical projection of $G$ on $G/H$ may not be a lattice of $G/H$. It suffices to take $P=Span \{e_{1},e_{2}+\alpha e_{3}\}$ where  $\{e_{1},e_{2},e_{3}\}$ is a set of generators for $\mathbb{Z}^{3}$ with  $\alpha$ irrational .

It is clear that \textbf {$L:=Ker (p\circ j)$ } is a subgroup of $H_{n}$ which contains $L_{1}$ as a finite index subgroup. Hence $L$ is a lattice of $H_{n}$ and we have the following split exact  sequence of groups:
\begin{center}
$1\leftharpoonup L\rightarrow \Gamma \rightarrow \mathbb{Z}\rightarrow 1$ 

\end{center}

\textbf{Theorem 3.} Every lattice $\Gamma$ of the Lie group $F_{n}$ is isomorphic to a group $\Gamma_{g(1)}$  semi-direct product of the normal Abelian group $L=Ker (p\circ j)$ by the group $ \mathbb{Z}$, via the action of the last given by a unipotent matrix $g(1)=I_{n}+\displaystyle\sum_{k> l} a_{kl}E_{k,l}$ in an ordered basis  $(e_{1},e_{2},...e_{n})$ of $L$.

Since the map $\mathbb{Z}\rightarrow Aut(L)$, $ p\mapsto g(p)=g(1)^{p}$, is a homomorphism of groups, the lattices in Theorem 3 are only determined by $g(1)$.

\textbf{Example 1}. If  $ dim F_{n}=3 $ then $ g(1)=\begin{pmatrix}

\ 1 & 0 \\
\ a & 1 
\end{pmatrix}$.
Here the isomorphism classes of lattices of $F_{n}$ are parametrized by the $|a|$ , $a\in \mathbb{Z}$.


Let $\Gamma=L\times \mathbb{Z}$, $\Gamma'=L'\times \mathbb{Z}$ be two lattices in $F_{n}$ and $\phi:\Gamma\rightarrow \Gamma'$  an isomorphism of lattices . From  Theorem 3, it is clear that $\phi=\varphi\times id $ where $\varphi:L\rightarrow L'$ is an isomorphism of groups.

As $ (v,p)(w,q)=(v+g(1)^{p}w,p+q)$ and $(\varphi(v),p)(\varphi(w),q)=(\varphi(v)+g'(1)^{p}\varphi(w),p+q)$, we have, as $\varphi$ is an isomorphism of groups:
\begin{center}
$(\varphi(v)+g'(1)^{p}\varphi(w),p+q)=(\varphi(v)+ \varphi(g(1)^{p}w),p+q)$
\end{center}
\begin{flushleft}
and, in particular, $g'(1)\varphi(w)=\varphi(g(1)w)$,i.e.
\end{flushleft}
\begin{center}

$\varphi^{-1}g'(1)\varphi=g(1)$
\end{center}
Conversely, this condition determines an isomorphism between $\Gamma$ and $\Gamma'$. Moreover, this condition implies $ \varphi^{-1}(g'(1)-I_{L'})\varphi=g(1)-I_{L}$ and more generally,
\begin{center}
$(g'(1)-I_{L'})^{k}\circ\varphi=\varphi\circ(g(1)-I_{L})^{k}$
\end{center}
for any $k\in\mathbb{N}$.

Notice that if the lattices associated to $g(1)$ and $g'(1)$ are isomorphic, these matrices are conjugated in $GL(n,\mathbb{Z})$ ,via a triangular matrix, with $\epsilon_{i}\in \{1,-1 \}$ in the diagonal. 

Denote by $L_{k}$ the kernel of $(g(1)-I_{L})^{k}$ ; idem for $L'_{k}$. We have $\varphi(L_{k})=L'_{k}$. In conclusion it results:

\textbf{Proposition 1.} Two lattices $L\times_{g(1)}\mathbb{Z}$ and $L'\times_{g'(1)}\mathbb{Z}$ in $F_{n}$ are isomorphic if an only if , in adapted basis to the above associated flags in $L$ and $L'$, we have  $ g(1)=  \varphi^{-1}\circ g'(1)\circ\varphi $ with $\varphi$ in the semi-direct product $T^{-}_{n}(\mathbb{Z})\times Diag( \epsilon_{i})$ where $T^{-}_{n}(\mathbb{Z})$ is the group of invertible lower triangular matrices with $1$ in the diagonal and coefficients in $\mathbb{Z}$.

\textbf{Remark 4.} Denote by $\{\Gamma_{g}\}$ the set of lattices in $F_{n}$; from the above study we have:

1. $(g(1)-I_{n})^{k}\neq 0 $ for $k$ integer, $ 0\leq k\leq n-1$ or equivalently $ g(1)=Exp(d')$ with $d'$ a $n$-nilpotent rational matrix of maximal nilpotent index. 

2. Let $N(T^{-}_{n}(\mathbb{Z}))$ be the normalizer of $T^{-}_{n}(\mathbb{Z})$ in $GL(n,\mathbb{Z})$. The group $N(T^{-}_{n}(\mathbb{Z}))$ acts over the set $\{g(1)\}$ and the orbits of this action are in bijection with the set of lattices $\{\Gamma_{g}\}$.

3. Consider the correspondence 
\begin{center}
$\theta:\{\Gamma_{g}\}\rightarrow ({\mathbb N}^{*})^{n-1}$ , $\Gamma_{g}\mapsto (|a_{2,1}|, |a_{3,2}|,...,|a_{n,n-1}|)$
\end{center}
This correspondence is in fact a map. The group $N(T^{-}_{n}(\mathbb{Z}))$  is generated by the matrices $U_{ij}=I_{n}+E_{i,j}$ with $i>j $ and the matrices $ Diag (\epsilon_{1},\epsilon_{2},...,\epsilon_{n})$. Via the conjugate action 
\begin{center}
$Diag (\epsilon_{1},\epsilon_{2},...,\epsilon_{n})g(1)Diag (\epsilon_{1},\epsilon_{2},...,\epsilon_{n})$
\end{center}
the coefficient $a_{ij}$ becomes $\epsilon_{i}\epsilon_{j}a_{ij}$ and to
the set  $\{a_{2,1},a_{3,2},...,a_{n,n-1}\}$ corresponds the set $\{\epsilon_{1}\epsilon_{2}a_{2,1},\epsilon_{2}\epsilon_{3}a_{3,2},...,\epsilon_{n-1}\epsilon_{n}a_{n,n-1}\}$. Hence the $|a_{i,j}|$ are not changed by the action of the matrices $Diag (\epsilon_{1},\epsilon_{2},...,\epsilon_{n})$.

\textbf{Example 2}. If $n=3$, $ g(1)= 
\begin{pmatrix}
1 & 0 & 0 \\
a & 1 & 0 \\
c & b & 1
\end{pmatrix}$
with $a,b,c \in \mathbb{Z}$
and $\varphi=
\begin{pmatrix}
1 & 0 & 0 \\
u & 1 & 0 \\
0 & v & 1
\end{pmatrix}$
with $u,v \in \mathbb{Z}$, because $\begin{pmatrix}
1 & 0 & 0 \\
0 & 1 & 0 \\
w & 0 & 1
\end{pmatrix}$ ,
with $w \in \mathbb{Z}$, is an element of $Z(T^{-}_{3}(\mathbb{Z}))$. Then $\varphi^{-1}=
\begin{pmatrix}
1 & 0 & 0 \\
-u & 1 & 0 \\
uv & -v & 1
\end{pmatrix}$ , $g(1)\varphi=
\begin{pmatrix}
1 & 0 & 0 \\
u+a & 1 & 0 \\
c+bu & b+v & 1
\end{pmatrix}$, and $\varphi^{-1}g(1)\varphi=
\begin{pmatrix}
1 & 0 & 0 \\
a & 1 & 0 \\
c+bu-av & b & 1
\end{pmatrix}$.

A simple observation of the coefficient $c':= c+bu-av $ permits the choice of $u$ and $v$ so that $0\leq c'< gcd (a,b)\leq Min(a,b)$. Hence if $ gcd (a,b)=1$ there is only one lattice, whereas if $gcd (a,b)= a>0$ there are $a$ lattices.

As a consequence, two lattices corresponding to $c'_{1}$ and $c'_{2}$ are isomorphic if and only if $c'_{1}+ (bu-av)=c'_{2}$.


We can prove :

\textbf{Proposition 2.} For any point $p\in({\mathbb N}^{*})^{n-1}$ and $\theta$  as in Remark 4, the set $\theta^{-1}(p)$ is finite. More precisely:

Given $g(1)=\begin{pmatrix} 
\ 1 & 0 \\ 
(a_{ij}) & 1 
\end{pmatrix}$ as in Theorem 3, there exists $\varphi\in T^{-}_{n}(\mathbb{Z})$ such that $\varphi^{-1}g(1)\varphi=
\begin{pmatrix}
\ 1 & 0 \\
(a'_{ij}) & 1 
\end{pmatrix}$ with $ 0\leq a'_{i,j}<a'_{j+1,j}$ for $i>j+1$.

\textbf{Proof.} Denote , in the proof, by $g$ the matrix $g(1)$ and by $g(e_{i})$ the $i$-column of $g$. First of all, by conjugation, we can suppose the coefficients of the first subdiagonal of $g(1)=\begin{pmatrix} 
\ 1 & 0 \\ 
(a_{ij}) & 1 
\end{pmatrix}$ strictly positive.  Order the elements $a_{i,j}$ as follows:
\begin{center}
$a_{3,1},a_{4,2},...,a_{n,n-2},a_{4,1},a_{5,2},...,a_{n,n-3},a_{5,1},...,a_{n,1}$
\end{center}
A Euclidean division gives $ a_{3,1}=qa_{2,1}+a'_{3,1}$ with $a'_{3,1}\in[0,a_{2,1}-1]$ and the  first three columns becomes $ g(e_{1})=e_{1}+a_{2,1}e_{1}+(qa_{2,1}+a'_{3,1})e_{3}+...$ ,
$ g(e_{2})=e_{2}+a_{3,2}+...$  and $ g(e_{3})=e_{3}+...$.Take as a new basis $\{e'_{i};0\leq i\leq n\}$ with $e'_{2}:=e_{2}+qe_{3}$ and $e'_{i}:=e_{i}$ if $i\neq 2$. We have then
$g(e'_{1})=e'_{1}+a_{2,1}e'_{2}+a'_{3,1}e'_{3}+a_{4,1}e'_{4}+...$ ; $g(e'_{2})=e'_{2}+a_{3,2}e'_{3}+...$. After an iteration the Euclidean division for $a_{4,2}$ and $a_{3,2}$ etc,we can hence suppose that $a_{3,1}\in[0,a_{2,1}-1],...,a_{n,n-2}\in[0,a_{n-1,n-2}-1],...,a_{i,j}\in[0,a_{j+1,j}-1]$.

The next step concerns $a_{i+1,j+1}$ if $i<n$ or $a_{n-j+2,1}$ if $i=n$. In the first case, if we make the Euclidean division $a_{i+1,j+1}=qa_{j+2,j+1}+a'_{i+1,j+1}$, the $j+1$-column of $g$ can be writen
\begin{center} $g(e_{j+1})=e_{j+1}+a_{j+2,j+1}e_{j+2}+...a_{i,j+1}e_{i}+(qa_{j+2,j+1}+a'_{i+1,j+1})e_{i+1}+...$
\end{center}

Consider the new basis $B'=\{e'_{i} ; 1\leq i\leq n\}$ with $e'_{j+2}:=e_{j+2}+qe_{i+1}$ and $e'_{i}:=e_{i}$ for the other indices. In the basis $B'$ the only coefficient changed in the $j+1$-column $g(e'_{j+1})$ is in the $i+1$ line where $a'_{i+1,j+1}$ replaces $a_{i+1,j+1}$. In the other columns the coefficients which change are after $a_{i+1,j+1}$ in the order adopted for the matrix terms.

The second case is treated as for $a_{3,1}$ and the proof is over.

 \hfill $\square$\\
 
 The lattices in $F_{n}$ with $a_{i+2,i+1}=1$ for every $i$, are all isomorphic and a basis change, as in the proof, reduces to the case where all $a_{i,j}$ with $i>j+1$ are zero : this lattice is $N(\mathbb{Z})$.

For this lattice $\Gamma_{0}$, we have $C_{i}\Gamma_{0}=C^{n-i}\Gamma_{0}$ for $0<i<n$, as one can easily verify using the presentation of $\Gamma_{0}$ :
\begin{center}
$\Gamma_{0}=<y_{1},...,y_{n},z>$  with $y_{i}y_{j}=y_{j}y_{i}$  for every $i,j$ 

$zy_{i}=y_{i}y_{i+1}z$  for $i<n$ and  $zy_{n}=y_{n}z$.

\end{center} 

For a general lattice $\Gamma$, the group $C^{n-i}\Gamma$ is contained in $C_{i}\Gamma$ for $0<i<n$ and the factor groups $C_{i}\Gamma/C^{n-i}\Gamma$ are finite Abelian groups. Their cardinalities can be computed in terms of the coefficients $a_{j+1,j}$. These assertions are obtained using the following presentation of $\Gamma$ as semidirect product:

\begin{center}
$\Gamma=<y_{1},...,y_{n},z>$  with $y_{i}y_{j}=y_{j}y_{i}$  for every $i,j$ ;  $zy_{n}=y_{n}z$ and for $i<n$,   $zy_{i}=y_{i}z_{i}z$ , where $z_{i}=y_{i+1}^{a_{i+1,i}}y_{i+2}^{a_{i+2,i}}...y_{n}^{a_{n,i}}$

\end{center}

Notice that if $C_{i}\Gamma=C^{n-i}\Gamma$ for all $i$, the lattice $\Gamma$ is isomorphic to $\Gamma_{0}$. In fact if $C^{1}\Gamma=C_{n-1}\Gamma$, we already have $\Gamma\cong\Gamma_{0}$. Moreover, two nonisomorphic lattices corresponding to matrices with the same first subdiagonal, may have isomorphic factor groups. 

\textbf{Example 3}. Let $n=3$ and consider the lattice $\Gamma$ with $g(1)$ as in \textbf{Example 2} and $0\leq c< gcd (a,b)\leq Min(a,b)$. The factor group $C_{1}\Gamma/C^{2}\Gamma$ is isomorphic to $C_{ab}$ independently of $c$. The other factor group
$C_{2}\Gamma/C^{1}\Gamma$ is the quotient of $\mathbb{Z}^{2}$ by the lattice generated by $(a,c)$ and $(0,b)$. Moreover, it is easy to see that if $\delta=gcd(a,b,c)$, this factor group is $C_{\delta}\times C_{ab\delta^{-1}}$. This shows that two nonisomorphic lattices may have isomorphic factor groups. For instance take $a=6, b=9, c_{1}=1$ for $\Gamma_{1}$ and $a=6, b=9, c_{2}=2$ for $\Gamma_{2}$; since $\delta=1$ for both of them, $C_{2}\Gamma/C^{1}\Gamma$ is cyclic of order $54$ for both lattices.

\section { Some Geometry of the corresponding affine compact nilmanifolds}

The Lie groups considered in section 2 are symplectic Lie groups (see [13]) and hence their left invariant symplectic form $ \Omega^{+}$  defines a left invariant affine structure $\nabla$ on the group by the formula (1). Since these groups are unimodular, this connection is geodesically complete ([6]). Consequently the quotient manifold $M=\Gamma\backslash G $, where $\Gamma$ is a lattice in $ G $, is a compact manifold endowed with  a symplectic form and a complete flat and torsion free linear connection $\nabla$. In fact, as these Lie groups are quadratic ( i.e. are endowed with a bi-invariant semi-Riemannian metric ) and the connection $\nabla$ on $G$ is the Levi-Civita connection of a flat left invariant semi-Riemannian metric ([7]), then $M$ carries also a flat semi-Riemannian metric.   

In the other hand, a direct calculation shows that the $2n$-dimensional Lie algebra ${\cal L}_{2n-1}(\mathbb{R})$, with bracket $ [e_{0},e_{i}]=e_{i+1}$ for $i=1,2,...,2n-2$ , has a scalar non degenerate cocycle given by $\Omega=\sum_{i=0}^{n-1}(-1)^{i}e_{i}^{*}\wedge e_{2n-i-1}^{*}$. Hence the groups $F_{2n-1}$ are symplectic Lie groups. Moreover as  any Lie algebra ${\cal L}_{k}$ has a non singular derivation $D$ its Lie bracket is compatible with the left symmetric product $ab=L_{a}b$ where $L_{a}=D^{-1}ad_{a}D$. Recall also that if $F_{n}$ is odd dimensional its carries a left invariant contact form ([13]). 

In short, in the following, we deal with \textbf{ affine compact and geodesically complete nilmanifolds}.The Euler characteristic of a such manifold is always zero (Kostant and Sullivan).

Let us now study \textbf{some special diffeomorphisms of the manifolds $M=\Gamma\backslash G $}. As the groups are nilpotent and simply connected, an automorphism of the Lie group $ G $, preserving a lattice  $\Gamma$, is fully determined by its restriction to $\Gamma$ (see [16]). Obviously any automorphism of $G$ stabilizing $\Gamma$, determines a diffeomorphism of $M$.

Consider the Heisenberg Lie group $G:=H_{1}(\mathbb{K})$ where $\mathbb{K}$ is $\mathbb{C}$, $\mathbb{R}\times \mathbb{R}$ or $\mathbb{D}$. The map $\phi:G\rightarrow G$  where the image by $\phi:=\phi_{\alpha,\beta,\gamma}$ of the matrix  
\begin{center}
$\begin{pmatrix}

1 & x & z \\   
0 & 1 & y \\
0 & 0 & 1

\end{pmatrix}$   is the matrix 
$\begin{pmatrix}

1 & \alpha x & \gamma z \\   
0 & 1 & \beta y \\
0 & 0 & 1

\end{pmatrix}$ 
\end{center}
\begin{flushleft}is an automorphism of $G$ if and only if $ \gamma=\alpha\beta$ with $\alpha,\beta,\gamma$
invertibles in $\mathbb{K}$ 
\end{flushleft}

Let $K$ a quadratic number field and $\Gamma$ the lattice of $G$ where $x,y,z$ are integers of $K$. The map $\phi$ will give an automorphism of $G$, preserving $\Gamma $, if and only if $\alpha,\beta,\gamma$ are invertible elements of the ring ${\cal O}_{K}$ of integers in $K$. Recall that if $K$ is complex, the group of units of ${\cal O}_{K}$ is finite: $C_{2}, C_{4}$ or $C_{6}$. When $K$ is real, the group of units is infinite and  direct product of $C_{2}$ by $\mathbb{Z}$, cyclic group generated by a so-called fundamental unit $\epsilon$. This means that any unit of ${\cal O}_{K}$ can be written as $\epsilon^{n}$ or $-\epsilon^{n}$with $n$ in $\mathbb{Z}$. For example, for $K$ the quadratic fields $\mathbb{Q}(\sqrt{k})$ with $k=2,3,5$ the corresponding $\epsilon$ is respectively $1+\sqrt{2}, 2+\sqrt{3}$ and $\frac{1+\sqrt{5}}{2}$.

Let $f$ be an arbitrary automorphism  of the Lie group $G$. Obviously $f|Z(G)$ is an automorphism of its center $Z(G)\equiv \mathbb{R}^{2}$ and its characteristic polynomial  has two roots. For $f=\phi_{\alpha,\beta,\gamma}$ the map $f|Z(G)$ corresponds to the map $ K\otimes_{\mathbb{Q}}\mathbb{R}\rightarrow K\otimes_{\mathbb{Q}}\mathbb{R}$, $z\otimes\lambda\mapsto \gamma z\otimes\lambda $ and the eigenvalues are the values of $\gamma$ for the two embeddings of $K$ in $\mathbb{R}$ that is  $\gamma(\sqrt{d})$ and $\gamma(-\sqrt{d})$.

In the other hand, $f$ induces an automorphism $\overline{f}$ of the 4-dimensional Abelian Lie group $G/Z(G)$ and hence has four characteric roots. In the case $f=\phi_{\alpha,\beta,\gamma}$, the map $\overline{f}:K^{2}\otimes_{\mathbb{Q}}\mathbb{R}\rightarrow K^{2}\otimes_{\mathbb{Q}}\mathbb{R}$ is the extension to $\mathbb{R}$ of the correspondence $(x,y)\mapsto(\alpha x,\beta y)$ and the eigenvalues of $\overline{f}$ are the images of $\alpha$ and $\beta$ by the real embeddings of $K$.

If $\alpha=\epsilon^{n}$ and $\beta=\epsilon^{m}$ then $f=\phi_{\alpha,\beta,\gamma}$, six eigenvalues. Moreover if $nm(n+m)\neq 0$, three of these eigenvalues have modulus $>1$ and the others modulus $<1$ and consequently $f$ induces an \textbf{Anosov diffeomorphism} of the space of right cosets $\Gamma\backslash G $.

In summary , if $U$ denotes the group of units of the ring ${\cal O}_{-d}$ and $Aut_{\Gamma}G$ is the group of automorphisms of $G$ preserving $\Gamma$ we have:

\textbf{Lemma 3.} The map $(\alpha,\beta)\mapsto \phi_{\alpha,\beta,\alpha\beta} $ is an injective homomorphism of groups from $U\times U$ into $Aut_\Gamma G$. In the complex case , $U$ is finite and  the eigenvalues of $\phi_{\alpha,\beta,\alpha\beta} $ have modulus 1. In the real case , $U$ is infinite and eigenvalues have inverse modules by pairs.

\textbf{Automorphisms of the Lie group $G:= T^{*}H_{1}(\mathbb{R})$}. In fact, it is easy to describe all the automorphisms of the Lie group $G$. We will only look for the automorphisms of  the lattice $\Gamma:=\Gamma_{1,1,1}$ using its presentation
\begin{center}
$y_{2}y_{3}=y_{3}y_{2}z_{1}$ , $y_{3}y_{1}=y_{1}y_{3}z_{2}$ , $y_{1}y_{2}=y_{2}y_{1}z_{3}$
\end{center}  
given in section 2. Consider three elements $y'_{1},y'_{2},y'_{3}$ in $\Gamma$ such that their canonical images in $\Gamma/Z\Gamma$ generate a lattice of $G/ZG$. This means that we can write 
\begin{center}
\begin{equation}
y'_{1}=y_{1}^{a_{1}}y_{2}^{a_{2}}y_{3}^{a_{3}}z'_{1} ,
y'_{2}=y_{1}^{b_{1}}y_{2}^{b_{2}}y_{3}^{b_{3}}z'_{2} ,
y'_{3}=y_{1}^{c_{1}}y_{2}^{c_{2}}y_{3}^{c_{3}}z'_{3} 
\end{equation}
\end{center}  
where $a_{i},b_{j},c_{k}$ are integers. The elements  $y'_{1},y'_{2},y'_{3}$  generate $\Gamma$ if and only if the matrix 
\begin{center}
$\begin{pmatrix}

a_{1} & b_{1} & c_{1} \\   
a_{2} & b_{2} & c_{2} \\
a_{3} & b_{3} & c_{3}

\end{pmatrix}$
\end{center}
is in $GL(3,\mathbb{Z})$ and the $z'_{i}$ are arbitrary elements of $Z(\Gamma)$. The commutators $<y'_{2},y'_{3}>, <y'_{3},y'_{1}> $ and $<y'_{1},y'_{2}> $ are in $Z(\Gamma)$ and generate it. Thus there is an unique automorphism of $\Gamma$ such that the image of $y_{i}$ is $y'_{i}$ for $i=1,2,3 $.

Hence, using the presentation of $\Gamma:=\Gamma_{1,1,1}$ given in the proof of Theorem 2, we have :

\textbf{Lemma 4.} For $y'_{i}$, with  $i=1,2,3$, in $\Gamma:=\Gamma_{1,1,1}$ as in (8), $z'_{i} \in Z(\Gamma)$ and the above matrix in $GL(3,\mathbb{Z})$ , there exists a unique automorphism $\phi$ of $\Gamma$ such that $\phi(y_{i})=y'_{i}$.

Let $f$ be an automorphism of $\Gamma:=\Gamma_{1,1,1}$. Denote by $f|Z(\Gamma)$ its restriction to the center of $\Gamma$ and by $\overline{f}$ the canonical map induced by $f$ on $ \Gamma/Z(\Gamma)$. Since $Z(\Gamma)$ and $ \Gamma/Z(\Gamma)$ are free $\mathbb{Z}$- modules of rank $3$ , we have associated to $f$ two invertible matrices $ A$ and $B$ with entries in $\mathbb{Z}$. Taking basis $ (z_{1}, z_{2},z_{3})$ and $(\overline{y_{1}},\overline{y_{2}},\overline{y_{3}})$ of these groups such that $ <y_{2},y_{3}>=z_{1} ,     ,    $, from the relations between $y_{i}$ and $z_{j}$ we have $A=\Lambda^{2}B$. As $det B =1$ or $det B=-1$, it follows $A= (det B) B^{-1}$. Hence, the characteristic polynomials of $B$ and $A$ are respectively given by 
\begin{center}
$ p_{B}(X)=X^{3}-(TrB)X^{2}+(Tr \Lambda^{2}B)X-detB $  

$q_{A}(X)=X^{3}-(Tr \Lambda^{2}B)X^{2}+((detB)^{-1}TrB)X-1$
\end{center}
  
Consequently we obtain :

\textbf{Lemma 5.} Let $f$ be an automorphism of the lattice $\Gamma:=\Gamma_{1,1,1}$ of the group $G:= T^{*}H_{1}(\mathbb{R})$. The six eigenvalues associated to the automorphisms $f|Z(\Gamma)$ and the induced one on $\Gamma/Z(\Gamma)$ , are algebraic integers. The product of the eigenvalues of $B$ is $1$ or $-1$ whereas the product of the eigenvalues of $A$ is always $1$ and up to the sign they are the inverses of the preceeding ones. If none have modulus unity, half have modulus bigger than one and the other modulus smaller than one. This happens if and only if no root of unity appears as eigenvalue.

\textbf{ Example 4}. Consider the matrix  of $GL(3,\mathbb{Z})$,
\begin{center}

$\begin{pmatrix}
1 & 5 & 2 \\   
2 & -1 & -1 \\
3 & 2 & 0
\end{pmatrix}$ 

\end{center}

\begin{flushleft}
Its characteristic polynomial is $ X^{3}-15X-1$ which has no root of modulus $1$. Hence the automorphism of $G:= T^{*}H_{1}(\mathbb{R})$ associated to this matrix determines an \textbf { Anosov diffeomorphism } of the manifold $\Gamma\backslash G$. 
\end{flushleft}

\textbf{Automorphisms of the Lie group $G:=F_{n}$}.

Let $\Gamma_{0}:=\mathbb{Z}^{n}\times \mathbb Z $ be the lattice of $G:=F_{n}$ , described in the subsection 3.2. An automorphism $ f $ of $\Gamma_{0}$ leaves invariant its central descending sequence. In terms of the generators $ y_{1}, y_{2},...y_{n},z $ this means that the Abelian subgroup $C^{1}\Gamma_{0}=<y_{2},...,y_{n}>$ and the flag $ <y_{i},y_{i+1},...,y_{n}>$ are invariant by $f$. Hence $f$ acts on the subgroups of $\Gamma_{0}$ containing $C^{1}\Gamma_{0}$. But all Abelian such subgroups are contained in
${\cal{M}}:= < y_{1},...,y_{n}>$ ; consequently $\cal{M}$ is invariant by $f$. So the restriction of $f$ to $\cal{M}$, in the $\mathbb{Z}$-basis $(y_{1},...,y_{n})$, is given by a lower triangular 
matrix with integers entries and diagonal terms $1$ or $-1$. Also, $f$ induces an automorphism of $\Gamma_{0}/\cal{M}$ and hence the image of $z$ by $f$ is $zm$ or $z^{-1}m$ where $m$ is any element of $\cal{M}$.

Notice , in particular, that the eigenvalues of $f$ are $1$ or $-1$. Putting $f(y_{i})=y_{i}^{\epsilon_{i}}$ with $\epsilon_{i}\in\{1,-1\}$ and $f(z)=zm$ and applying $f$ to the relation $zy_{i}=y_{i}y_{i+1}z$ we get $zmy_{i}^{\epsilon_{i}}=y_{i}^{\epsilon_{i}}y_{i+1}^{\epsilon_{i+1}}zm $. As $my_{i}=y_{i}m $ then $zy_{i}^{\epsilon_{i}}=y_{i}^{\epsilon_{i}}y_{i^+1}^{\epsilon_{i+1}}z$ i.e. $\epsilon_{i}=1$ implies $\epsilon_{i+1}=1$. In a similar way $\epsilon_{i}=-1$ implies $\epsilon_{i+1}=-1$.

\textbf{Toral Affine Symplectic actions }.

By definition, the  action $L^{G}$ on a symplectic Lie group $(G,\Omega^{+})$ given by $L^{G}_{\sigma}: \tau\mapsto \sigma\tau$, is symplectic. So, the right invariant vector fields $x^{-}$ with  $x\in {\cal G}$ on $G$ are symplectic ( or locally Hamiltonian ). If $G$ is simply connected, these vector fields are Hamiltonian , that is, the action $L^{G}$ is Hamiltonian. In fact , in this case, a direct calculation ( see [10])shows that the map $Q
:exp x\mapsto \sum_{k=1}^{\infty} \frac{1}{k!}( ad^{*}_{x})^{k-1}.\Omega(x,.)$ is an equivariant moment map for $L^{G}$ and the action $\rho_{G}$ of $G$ on ${\cal G}^{*}$ given by $\sigma\mu:=Q(\sigma)+Ad^{*}(\sigma)(\mu)$. In particular we have :
\begin{center}
$i(x^{-})\Omega^{+}=-dQ_{x}$ , where $Q_{x}(expy)=<Q(expy),x>$ with $x,y \in {\cal G}$.
\end{center}
As $Q(\sigma\tau)= Q(\sigma)+Ad_{\sigma}^{*}(Q(\tau))$ for any $\sigma,\tau$ in $G$ we will say ([10]) that $Q$ is moment cocycle. Moreover since , $\sigma=exp x\mapsto (Q(\sigma),Ad^{*}(\sigma))$ is a representation of $G$ by affine transformations of $ {\cal G}^{*}$ and  $Q$ is a diffeomorphism , $Q$ is a  developing map of the affine structure on $M:=\Gamma\backslash G$ with $\Gamma$ any lattice of $G$ and $\nabla$ is geodesically complete. Notice that $Q(\epsilon)=0$, where $\epsilon$ stands for the unit of $G$  , and $h:\gamma\mapsto h(\gamma)=(\gamma,Q(\gamma))$ is the holonomy representation of the affine structure of $M$. We have then:
\begin{center}
$h(\gamma)\circ Q = Q\circ \gamma$ for any $\gamma\in\Gamma$ 
\end{center} 
So, $M=\Gamma\backslash G$ can be identified to the manifold $h(\Gamma)\backslash{\cal G}^{*} $, with $\Gamma$ acting properly and freely on ${\cal G}^{*}$ via $h$. The tangent bundle $TM$ identifies with the fiber product 
\begin{center}
$E_{h}:=G\times_{h}{\cal G}^{*}=(G\times {\cal G}^{*})/\Pi_{1}(M)$
\end{center}
\begin{flushleft}
where $\Pi_{1}(M)=\Gamma $ acts diagonally by deck transformations on the factor $G$ and via $h$ on the factor ${\cal G}^{*}$.
\end{flushleft}
 The developing map $Q$ defines an isomorphism of vector bundles of $TM$ with $E_{h}$. The flat connection $\nabla$ on $TM$ arise from the representation $h$ in the standard way. For the structure of symplectic Lie groups the reader can refer to [10]. Recall that $M$ is endowed with a Lagrangian foliation (see [17]). 

Let $ (u_{1},u_{2},...,u_{n})$ a basis of $ {\cal G}$ and $p:G\rightarrow M$ the canonical projection. It is clear that $p$ is a symplectic and affine covering map. The global system of affine coordinates $Q_{u_{i}}(\sigma):=Q(\sigma)(u_{i})$ for $i=1,2,...n$ on $G$ determines a local system of affine coordinates $(x_{1},x_{2},...,x_{n})$ on $M$ given by $\Omega^{*}(X_{i},.)=dx_{i}$, where $ \Omega^{*}$ is the symplectic form on $M$ and $X_{i}$ is a local vector field on $M$ verifying $p_{*}(u_{i}^{-})=X_{i}$. Using the coordinates $x_{i}$, Boucetta and the second author have shown in [12] the following facts:

1) The Poisson bracket corresponding to $\Omega^{*}$ is polynomial of degree 1.

2) The symplectic form $\Omega^{*}$ is polynomial of degree at most $n-1$ where $n$ is the dimension of $M$.

3) The volume form $(\Omega^{*})^{n}/n!$ is parallel relative to $\nabla$.

4) Any solution of the classical Yang-Baxter equation $ {\cal G}$ gives arise to determines a polynomial Poisson tensor on $M$. 

In the following $G$ is nilpotent and simply connected and $\Gamma$ is a lattice in $G$. Notice that the action $\Psi$ of $Z(G)$ on $M$, deduced of $L^{Z(G)}$, is affine and symplectic. As the stabilizer of any point in $M$, under $\Psi$, is $Z(G)\cap\Gamma $, the torus $ T^{*}:=Z(G)/Z(G)\cap \Gamma $ acts freely on $M$. Consequently, using a Corollary of the slice Theorem of Koszul (see [18], page 180 ),  $M/T^{*}$ is a manifold and $M$ is a principal $T^{*}$ -bundle over $M/T^{*}$. Obviously the canonical projection of $M$ on $M/T^{*}$ is an affine map.   
 
In the other hand any closed subgroup $H$ of $G$ acts naturally, via the action $L^{H}$, by symplectomorphisms of $G$ and  the map $Q_{H}=i'\circ Q$ , with $i'$ the transpose of the canonical injection $i$ of the Lie algebra of $H$ in ${\cal G}$, is a $\rho_{H}$ equivariant invariant moment map relative to $L^{H}$ . Notice that any value of $Q_{H}$ is a regular value.


Recall now that the bracket of two symplectic vector fields is Hamiltonian and suppose $G$ non Abelian. Consequently if we take $H:=Z(G)\bigcap D(G)$ , the action $L_{T'}$ of the torus $T':= H/H\bigcap Z(\Gamma)$ on $M$ becomes Hamiltonian. More precisely if $X,Y$ are symplectic vector fields on $M$ then $[X,Y]$ is Hamiltonian with Hamiltonian function $ \Omega^{*}(Y,X)$. Hence we have (see [18]):

\textbf{Lemma 6.} The action of the torus $T':= H/H\bigcap Z(\Gamma)$ on $M=\Gamma\backslash G$ induced by $L^{H}$ is Hamiltonian if $G$ is nilpotent and non Abelian. In particular if $G=F_{n}$ there is a Hamiltonian action of a circle on $M$  and, following Duistermaat-Heckman, the stationary-phase approximation is in this case exact.

The determination of the Duistermaat-Heckman function in the Lemma could be interesting ([18],page 70).


Denote by $Q_{T'}$ the corresponding moment map in the Lemma 6. According to the convexity theorem of Atiyah and Guillemin-Sternberg ([8], page 59 ), the image of $Q_{T'}$ is the convex hull of the image of the $T'$ fixed point set of $M$. Moreover, in the case where dim $Z(G)\bigcap D(G)$ is the half of dim $ M $ the image of the moment map is  necessarily a Delzant polytope ([19],page 75)) because the action $\Psi_{T}$ is effective. This is the case for example if $G:= T^{*}H_{1}(\mathbb{R}) $.

\textbf{Theorem 4.} Let $(G,\Omega^{+})$ be a symplectic nilpotent non Abelian and simply connected Lie group and let $\tau\in Z(G)\bigcap D(G)\bigcap\Gamma$ with $\tau\neq\epsilon$, $\tau=expz$ and $H:=exp(\mathbb{R}z)$. Then the action $L_{T}$ of the cercle $T:=H/exp(\mathbb{Z}z )$ on $M$ given by $\widehat{\tau_{1}}.\overline{\sigma}:=\overline{\tau_{1}\sigma}$, deduced of $L^{H}$, is Hamiltonian and the corresponding moment map $Q_{T}$ is $\rho_{T}$ - equivariant. Moreover 
$N:= T\backslash M$ is an affine manifold, $T$ acts on $M$ by (local) translations and $M$ is an affine $T$-principal bundle over $N$. Finally, the origin of $(\mathbb{R}z)^{*}$ is a regular value of $Q_{T}$, the set level $Q_{T}^{-1}(o)$ is an affine hypersurface of $M$ and   
$R:=T\backslash Q_{T}^{-1}(o)$ is a reduced (affine) symplectic nilmanifold of $M$ by means of $T$. In fact we have the following canonical sequences of affine manifolds and affine maps :  

\begin{center}

$Q_{T}^{-1}(0)\leftharpoonup M\rightarrow (Lie(T))^{*}$

\end{center}

\begin{center}

$T\leftharpoonup Q_{T}^{-1}(0)\rightarrow R $ 

\end{center}

\textbf{Proof.} First of all notice that $H$ is a central subgroup, $H\bigcap \Gamma$ is a lattice in $H$ and $T$ is a circle. As above, the action of $T$ on $M$ is symplectic and Hamiltonian. Moreover, as $ p $ is a local symplectomorphism, if $ z^{*} $ is the symplectic vector field on $M$, we have $p_{*}(z^{-})= z^{*}$. Consequently the regular value $0$ of $Q_{H}$ is also a regular value of $Q_{T}$. Hence $Q_{T}^{-1}(0)$ is a hypersurface of $M$. A simple verification shows that the action of $T$ on $M$ is free. Using a $T$-invariant Riemannian metric on $M$ and the slice theorem of Koszul, it follows that $N:= T\backslash M$ is a manifold and $M$ is a $T$-principal bundle over $N$. The fact that $H$ acts on $G$ by translations ( see [10] ) implies that $T$ acts (in local coordinates ) on $M$ in a similar way. As a consequence, the canonical maps involved in the $T$-principal bundle are affine maps. Finally,the moment map $Q_{T}$ is $\rho_{T}$-equivariant and the stabilizer group $T_{0}$ , by the action $\rho_{H}$, acts on $Q_{T}^{-1}(0)$ freely and properly. Hence, using the procedure described by Marsden and Weinstein in [20], there exists a unique symplectic form $\Omega^{*}_{0}$ on $R$ with the property $\pi_{0}^{*} \Omega^{*}_{0}=i^{*}_{0}    \Omega^{*}$ where $\pi_{0}:Q_{T}^{-1}(0)\rightarrow R$ is the canonical projection and $i_{0}:Q_{T}^{-1}(0)\rightarrow M$ is the canonical inclusion.

\hfill $\square$\\

\textbf {Fundamental group of the reduced manifold.} Given a $2n$-dimensional compact symplectic nilmanifold $( \Gamma\backslash G,\Omega^{*})$, it is natural to ask if we can, in general, easily describe the fundamental group of a $2n-2$  symplectic reduced manifold $N$ of $M$ obtained , as shown by Theorem 4, in terms of $\Gamma $. The response to this question is in general negative, as shown the following analysis.

Recall some facts about the structure of simply connected nilpotent symplectic Lie groups. Denote by $H^{\bot}$ the connected component of the unit in $Q_{H}^{-1}(0)$. Then $H^{\bot}$ is a closed subgroup of $G$ which is $ \Omega^{*}$-orthogonal to $H$ and $K:=H\backslash H^{\bot}$ is a reduced symplectic Lie group. Moreover we have the exact canonical sequence of affine Lie groups
\begin{center}
$ \{\epsilon\}\leftharpoonup H\leftharpoonup H^{\bot}\rightarrow H\backslash H^{\bot}\rightarrow \{\epsilon\}$
\end{center}
and the canonical affine principal fiber bundle

\begin{center}

$ \{\epsilon\}\leftharpoonup H^{\bot}\leftharpoonup G\rightarrow  H^{\bot}\backslash G$

\end{center}
\begin{flushleft}
where the projection can be identified with the moment map $Q_{H}$ (see [10]).
\end{flushleft}

Let $\Gamma$ be a lattice in $G$. If $H^{\bot}$ admits lattices and the behaviour of $\Gamma$ relatively to $H$ and $H^{\bot}$ is nice i.e. $H\bigcap \Gamma$ and $H^{\bot}\bigcap\Gamma$ are lattices in $H$ and $H^{\bot}$, we obtain a compact symplectic reduced manifold $R'$ of $M$, with fundamental group $\Gamma':=(\Gamma\bigcap H)\backslash \Gamma\bigcap H^{\bot}$, by means of the compact affine manifold  $(H^{\bot}\bigcap\Gamma)\backslash H^{\bot}$ or the symplectic Lie group $H\backslash H^{\bot} $.

Notice that the manifolds involved in the last sequences are covering spaces of the manifolds involved in the sequences described in the Theorem 4.   

\textbf {Example 5.} Let $G$ be the Lie group given by the manifold $\mathbb{R}^{6}$ endowed with the product.
\begin{center}
$ (x_{i},y_{j})(x'_{i},y'_{j})= (x_{i}+x'_{i},y_{j}+y'_{j}+x_{k}x'_{l})$
\end{center}
where $\{j,k,l\}$ is a cyclic permutation of $\{1,2,3\}$. The 1-parameter subgroups of $G$ correspond to the maps
\begin{center}
$\varphi:t\rightarrow (tx_{i},ty_{j}+ \frac{t(t-1)}{2}x_{k}x_{l})$
\end{center}
and then

\begin{center}
$\frac{d\varphi}{dt}(0)=(x_{1},x_{2},x_{3},y_{1}-\frac{x_{2}x_{3}}{2}, y_{2}-\frac{x_{3}x_{1}}{2},y_{3}-\frac{x_{1}x_{2}}{2})$
\end{center}
and 
$\varphi(1)= (x_{1},x_{2},x_{3},y_{1},y_{2},y_{3})$.

Consequently the exponential map $ exp: {\cal G}\rightarrow G$ is given by 
\begin{center}
$(a_{1},a_{2},a_{3},b_{1},b_{2},b_{3})\rightarrow (a_{1},a_{2},a_{3},b_{1}+\frac{a_{2}a_{3}}{2},b_{2}+\frac{a_{3}a_{1}}{2},b_{3}+b_{1}+\frac{a_{1}a_{2}}{2})$

\end{center}
\begin{flushleft}
whereas the logarithm map is described by
\end{flushleft}
\begin{center}
$ln: G\rightarrow {\cal G}, (x_{i},y_{j})\longmapsto (x_{i},y_{j}-\frac{y_{k}x_{l}}{2})$

\end{center}

Let us look for the scalar 2-cocycle $\Omega$ on ${\cal G}$. As $D({\cal G})=Z({\cal G})$ then $\Omega$ is totally isotropic on $D({\cal G})$. Using the natural basis $\{e_{i},f_{j}\}$ of ${\cal G}$ we have the condition $\Sigma^{i=3}_{i=1}\Omega(e_{i},f_{i})=0$. Up to a 2-coboundary, the matrix of $\Omega$ can be written 
\begin{center}
$\begin{pmatrix}
0& -B^{t} \\   
B & -0 \\
\end{pmatrix}$ 
\end{center}
where $B$ is a square matrix of size 3 and trace 0, invertible if $\Omega$ is non degenerate.

Take for $\Gamma$ the lattice of $G$ with $x$ and $y$ integers. Any 1-dimensional central subgroup $H$ of $G$ which meets $\Gamma$ contains an element $z=(0,0,0,y_{1},y_{2},y_{3})$ where the $y_{i}$ are coprime integers. Using an automorphism of $G$ leaving $\Gamma
$ fixed, $H$ becomes $H=\{(0,0,0,0,0,t), t\in \mathbb{R}\}$. Putting $B=(b_{ij})$ we have then, $I=\mathbb{R}f_{3}$ and $I^{\bot}=V\oplus Z({\cal G})$ with
\begin{center}
$V=\{(a_{1},a_{2},a_{3}); a_{1}b_{31}+a_{2}b_{32}+a_{3}b_{33}=0\}$.

\end{center} 
\begin{flushleft}

Hence $ exp( I^{\bot})=\{(x_{1},x_{2},x_{3},y_{1},y_{2},y_{3}) ; \Sigma^{i=3}_{i=1}x_{i}b_{3i}=0\}$
\end{flushleft}
We want to describe $\Gamma':=\Gamma\bigcap exp( I^{\bot})$. If the $b_{3,i}$ are $\mathbb{Q}$-linearly independent, the $x_{i}$ are zero and $\Gamma'=Z(\Gamma)$. 

If the $\mathbb{Q}$-vector space  $W$ generated by the set $\{b_{3i}; i=1,2,3 \}$ has dimension 2, we have $(x_{1},x_{2},x_{3})= d(t_{1},t_{2},t_{3})$ with $d$ and $t_{i}$ in $\mathbb{Q}$. So , $\Gamma'=\{(nu_{1},nu_{2},nu_{3},v_{1},v_{2},v_{3}) ; n,v_{i}\in \mathbb{Z}\}$ where $(u_{1},u_{2},u_{3})$ is an unimodular element of $\mathbb{Z}^{3}$. Finally, if $W$ is 1-dimensional, the relation $ \Sigma^{i=3}_{i=1}x_{i}b_{3i}=0$ tells that there exist a non null linear form  $f$ on $ \mathbb{Z}^{3}$  such that 
\begin{center}
$ \Gamma'=\{(x_{1},x_{2},x_{3},y_{1},y_{2},y_{3}) ; x_{i},y_{j} \in \mathbb{Z}, f(x_{1},x_{1},x_{1})=0\}$
\end{center}

\textbf {In short $ \Gamma'$ is a lattice in the Lie group $H^{\bot}$ only in this last case.}

\textbf{Proposition 3}. Let $G=F_{2m-1}$ , $H=Z(G)$ and $\Gamma=L\times_{g(1)}\mathbb{Z}$ as in Proposition 1. Then $H^{\bot}$ do not depend of the left invariant symplectic form and  $\Gamma\bigcap H$ and $\Gamma\bigcap H^{\bot}$ are lattices in $H$ and $H^{\bot}$ respectively. In fact, $H^{\bot}$is the integral subgroup of $G$ associated to the unique Abelian maximal ideal of ${\cal G}$ and $\Gamma_{1}:=L/Ker (d-Id)$ is a lattice of the Abelian Lie group $K:=H\backslash H^{\bot}$. Also, the torus $R:=\Gamma_{1}\backslash K $ is a  reduced manifold of the symplectic nilmanifold $M:=\Gamma\backslash G$.

\textbf{Proof.} Let $\{e_{0},e_{1},...,e_{2m-1}\}$ be an  adapted basis of ${\cal G}$ in the sense of Vergne i.e. $ad_{e_{0}}(e_{i})=e_{i+1}$ for $0 < i< 2m-1$ . A two scalar cocycle $\omega $ on ${\cal G}$ verifies $\omega (e_{i},e_{j})=0$ whenever $i+j$ is even or bigger than $2m-1$. It follows that 
$ Z(\cal G)^{\bot}$ contains $V= Span \{e_{1},...,e_{2m-1}\}$ for every $\omega $. If $\omega $ is nondegenerate, $ Z(\cal G)^{\bot}$  has dimension $2m-1$, so is equal to $V$. This proves the first claim. The description of lattices of $F_{2m-1}$ in Proposition 1, implies the assertions concerning $\Gamma_{1}$.

In the other hand, the group $K$ is a symplectic Lie group, reduction of the symplectic Lie group $G$ by the Lie group $H$. Also, there exists a natural embedding, denoted in the following  by $i$, of the torus $N:=(\Gamma\bigcap H^{\bot})\backslash H^{\bot}$ into the manifold $M$ and the natural action of the cercle $C:=(\Gamma\bigcap H)\backslash H$  on $N$ determines the compact manifold $R:=C\backslash N $. Moreover, as this action is free (and proper) a Corollary of the slice Koszul Theorem (see [18]) produces  a $C$-principal fiber bundle 
\begin{center}

$ C\rightarrow N \stackrel{\pi}{\rightarrow} R   (*)$

\end{center}

Notice that the canonical action of $H$ on $H^{\bot}$ preserves the restriction of the symplectic form  $\omega^{+}$ to $H^{\bot}$ and this restriction induces a left invariant symplectic form $\omega_{K}$.
Also, the canonical action $\phi$ of $C$ on $K$ is symplectic and $J_{K}(\overline{exp x)}=\sum((ad_{x}^{*})^{k-1}.\omega(x,.))$ , with $x$ in the Lie algebra of $H^{\bot}$, is a moment map for $\phi$. Finally, we have directly ( or using the procedure Marsden-Weinstein ) a symplectic form $\Omega_{R}$ on $R$ such that $\pi^{*}(\Omega_{R})=i^{*}(\Omega)$. In conclusion, the torus $R$ is a symplectic manifold of dimension $2n-2$ reduction of the $2n$ symplectic nilmanifold $M$. 
Remark that that the objects and the maps involved in fibration (*) are affine maps. This also the case for 

\begin{center}

$N \stackrel{i}{\rightarrow}M\stackrel{J_{K}}{\rightarrow}J_{K}^{-1}(\nu)   (**)$

\end{center} 

\begin{flushleft}
where $\nu $ is a regular value of $J_{K}$.
\end{flushleft}
Hence we can say, as in [10], that the symplectic manifold $M$ is a	 symplectic double extension of the symplectic manifold $R$.\hfill $\square$\\

\section {Symplectic Heisenberg Lie groups}

In this section $\mathbb{A}$ is a local associative and commutative $l$- dimensional real algebra and ${\cal H}_{k}(\mathbb{A})$ is the corresponding $l(2k+1)$-dimensional Heisenberg Lie algebra defined by $\mathbb{A}$. We will give a necessary and sufficient condition for the existence of a scalar nondegenerate $2$-cocycle on ${\cal H}_{1}(\mathbb{A})$. More precisely, if we decompose $\mathbb{A}$ as $\mathbb{A}=\mathbb{R}\oplus\mathbb{N}$, where $\mathbb{N}$ is the maximal ideal of $\mathbb{A}$ and we denote by $\mathbb{S}$ the annihilator of $\mathbb{N}$ or socle of $\mathbb{A}$, we have :

\textbf{Theorem 5.} The real Lie algebra ${\cal H}_{1}(\mathbb{A})$ is a symplectic Lie algebra if and only if $\mathbb{A}$ is even dimensional and $d=dim  \mathbb{S}\leq 2$. Moreover for $k\geq 2$ the algebra ${\cal H}_{k}(\mathbb{A})$ is not symplectic.

\textbf{Proof.}The schema of the proof is as follows: if $d \geq 3$, we show that any scalar $2$-cocycle is degenerate ; if $ d=1$ then $\mathbb{A}$ is a Frobenius algebra and we describe a nondegenerate scalar $2$-cocycle. The case $ d=2$ requires a more detailed study.

It is clear that ${\cal H}_{1}(\mathbb{A})={\cal H}_{1}(\mathbb{R})\otimes \mathbb{A}= e\otimes \mathbb{A}\oplus f\otimes \mathbb{A}\oplus g\otimes \mathbb{A}$ with $[e,f]=g$. 
Suppose $d=1$ and let $\mu $ be a linear form on $\mathbb{A}$ non zero on $\mathbb{S}$ so that the  bilinear form $ \varphi (x,y)=\mu (xy) $ satisfying $ \varphi( xy,z)= \varphi (x,yz)$ is nondegenerate. A direct calculation shows that the formula

\begin{center}

$\omega (e\otimes a + f\otimes b+ g\otimes c,e\otimes a' + f\otimes b'+ g\otimes c')= \mu (ac'-a'c)+ \omega'(b,b')$

\end{center}
defines a non degenerate $2$-cocycle on ${\cal H}_{1}(\mathbb{A})$ if and only if $l$ is even and $\omega'(b,b')$ is a nondegenerate alternating  bilinear form.

If $d\geq 3$, a simple calculation shows that $\mathbb{N}e+\mathbb{N}f+\mathbb{A}g $ is orthogonal to 
$\mathbb{S}g$ relatively to any scalar $2$-cocycle of ${\cal H}_{1}(\mathbb{A})$. Consequently any $2$-cocycle  is degenerate.

Now, suppose $d=2$ and recall that an alternate bilinear form $\omega$ on ${\cal H}_{1}(\mathbb{A})$  is a $2$-cocycle if and only if $\mathbb{A}g $ is totally isotropic for $\omega$ and we have $\omega(ae,bg)=\omega(abe,g) $ ; $\omega(af,bg)=\omega(abf,g)$ for any a,b in $\mathbb{A}$.

Let $f_{1}, f_{2}$  two linear forms on  $\mathbb{A}$ such that their restrictions to $\mathbb{S}$ are independent and $\varphi_{1}, \varphi_{2}$ the associatives bilinear forms given by $ \varphi_{i}(a,b)=f_{i}(ab)$ for $i=1,2$.

Define a bilinear form  $\omega$ on ${\cal H}_{1}(\mathbb{A})$ putting: 1) $\omega(cg,c'g)=0$ , 2) $\omega(ae,bg)= \varphi_{1}(a,b)$ , 3) $\omega (af,bg)= \varphi_{2}(a,b) $, 4) The restriction of $\omega$ to $\mathbb{A}e\oplus \mathbb{A}f $ is any alternate bilinear form. In fact $\omega$ is a $2$-cocycle. Moreover we will  show that the linear map $\mathbb{A}g\rightarrow {\cal H}_{1}(\mathbb{A})^{*}$, given by $ cg\rightarrow \omega(cg,.)$ is injective. Consider $xg$ in the kernel of this map, we have $\varphi_{1}(a,x)=\varphi_{2}(a',x)=0$ for any $a, a'$ in $\mathbb{A}$ , so $f_{1}(ax)= f_{2}(a'x)=0$. Let us show that $x=0$. If $x\neq 0$, consider $I$ a minimal nonzero ideal of $\mathbb{A}$ contained in $\mathbb{A}x$. We have 
$\mathbb{N}I=\{0\}$, because $\mathbb{N}I$ is an ideal and $\mathbb{N}I=I$ would imply $I=\{0\}$. So, $I$ is contained in $\mathbb{S}$ and there exists $a\in \mathbb{A}$ such that $0\neq ax \in \mathbb{S}$. Then $f_{1}(ax)= f_{2}(ax)$ and as  the restrictions to $\mathbb{S}$ of $f_{1}$ and $f_{2}$ are independent,it follows $ax=0$, which is impossible.

Consequently, there exists a subspace $E$ of $\mathbb{A}e\oplus \mathbb{A}f$ in duality with $\mathbb{A}g$ by $\omega$. Denote by $F$ a supplementary subspace of $E$ in $\mathbb{A}e\oplus \mathbb{A}f$. According to condition 4) above,  we can choose the restriction of $\omega$ to $\mathbb{A}e\oplus \mathbb{A}f$ verifying $\omega(E,E)=\omega(E,F)=0$ and $\omega|F\times F$, any nondegenerate alternate form (notice that $ dim F=dim \mathbb{A}$ is even by hypothesis). This choice determines a nondegenerate $2$-cocycle of ${\cal H}_{1}(\mathbb{A})$.

Finally consider the case $k\geq 2$. Obviously, ${\cal H}_{k}(\mathbb{A})={\cal H}_{k}(\mathbb{R})\otimes \mathbb{A}$ and ${\cal H}_{k}(\mathbb{A})=\mathbb{A}e_{1}\oplus...\oplus\mathbb{A}e_{k}\oplus\mathbb{A}f_{1}\oplus...\oplus\mathbb{A}f_{k}\oplus\mathbb{A}g$ with $[e_{i},f_{i}]=g$. If $\omega$ is a $2$-cocycle of
${\cal H}_{k}(\mathbb{A})$, we have for any $a\in \mathbb{A}$:
\begin{center}
$\delta\omega(ae_{1},f_{1},e_{i})=0=\delta\omega(ae_{1},f_{1},f_{i})$ for $i \geq2$
\end{center}
and so $\omega(ag,e_{i})=\omega(ag,f_{i})=0$.

As $\delta\omega(ae_{2},f_{2},e_{1})=0=\delta\omega(ae_{2},f_{2},f_{1})$, we also have   $\omega(ag,e_{1})=\omega(ag,f_{1})=0$ and 
 $\delta\omega(ae_{1},f_{1},bg)=0$ implies $\omega(ag,bg)=0$. In summary, $\mathbb{A}g$ is  in the kernel of $\omega$. Hence $\omega$ is degenerate. \hfill $\square$\\
 
\textbf{Remark 5.} It is clear in the proof of Theorem 5 that ${\cal H}_{1}(\mathbb{A})$ is symplectic if $\mathbb{A}$ is any Frobenius algebra non necessary a local algebra.
 
\textbf{Example 6.} Let $\mathbb{A}:=\mathbb{R}[x,y]/(x^{3}=y^{2}=xy=0) $ , so this algebra can be decomposed as $\mathbb{A}=\mathbb{R}\oplus \mathbb{R}x \oplus \mathbb{R}x^{2}\oplus \mathbb{R}y$ and $\mathbb{S}=\mathbb{R}x^{2}\oplus \mathbb{R}y$ is its socle. Following step by step the proof of the Theorem 5 we can find a scalar nondegenerate cocycle on ${\cal H}_{1}(\mathbb{A})$.
 
\textbf{Example 7.} Let $\mathbb{A}=\mathbb{R}^{2k}\oplus \mathbb{C}^{h}$ and ${\cal K }:= \mathbb{Q}[X]/(P)$ where $P$ is a polynomial having $2k$ real roots and $2h$ non real roots . Then $ H_{1}({\cal O}_{\cal K})$ is a lattice in the real simply connected Lie group of Lie algebra ${\cal H}_{1}(\mathbb{A})$.

\textbf{Example 8.} Consider the real algebra $\mathbb{A}_{\cal Q}=\mathbb{R}\oplus V\oplus\mathbb{R}$, with $V$ a $2d$- dimensional vector space and product given by 
\begin{center}

$(\lambda,v,\mu)(\lambda',v',\mu')=(\lambda\lambda',\lambda v'+\lambda' v,\Phi_{\cal\ Q}(v,v')+\lambda\mu'+\lambda'\mu)$

\end{center}
\begin{flushleft}
where $\Phi_{\cal\ Q}$ is the bilinear symmetric form associated to a non degenerate quadratic form ${\cal Q}: V \rightarrow \mathbb{R}$. The subset $\{(0,v, \mu); v\in V, \mu\in \mathbb{R}\}$ is the maximal ideal of $\mathbb{A}_{\cal Q}$. It is also a Frobenius algebra via the linear form $f(\lambda,v,\mu)=\mu$. Take $ \overline{q}:\mathbb{Q}^{2d}\rightarrow \mathbb{Q}$ a rational quadratic form with $sign( \overline{q})=sign ({\cal Q})$. Then, the quadratic real extension $ {\overline{q}}$ is isometric to ${\cal Q}$. Finally, if we denote by $\mathbb{A}_{\overline{q}}$ the the corresponding algebra associated to $ \overline{q}$ as above, the real Lie algebras ${\cal H}_{1}(\mathbb{A}_{\overline{q}})\otimes_{\mathbb{Q}} \mathbb{R}$ and ${\cal H}_{1}(\mathbb{A}_{\cal Q})$ are isomorphic. This defines a lattice of the real simply connected  Lie group of which the Lie algebra is ${\cal H}_{1}(\mathbb{A}_{\cal Q})$.
\end{flushleft}

As two algebras $\mathbb{A}_{\overline{q}}$ and $\mathbb{A}_{\overline{q'}}$ are isomorphic if and only if $ \overline{q}$ and $ \overline{q'}$ are collinear, the commensurability classes of lattices in the group $H_{1}(\mathbb{A}_{\cal Q})$ are in bijection with the classes of collinear rational quadratic forms of the same index than ${\cal Q}$.

\section{Some symplectic solvmanifolds with  symplectic or semi-Riemannian affine structure} 

First of all recall that the left invariant connection defined by the formula (1) is not symplectic if $G$ is non Abelian. Nevertheless in [15], we have studied lattices in 4-dimensional Lie symplectic Lie groups and showed that these groups admit a complete left invariant  symplectic affine structure. Here we give a simpler proof of this result, in a more general framework.

\textbf{Theorem 6.} Let $(G,\omega^{+})$ be a symplectic Lie group having an Abelian  normal closed codimension one subgroup. Then $(G,\omega^{+})$ admits a complete left invariant symplectic affine structure. Consequently any homogeneous space $K\backslash G$ inherits a complete flat and torsion free symplectic connection. In particular this is the case if $K$ is a (uniform) lattice in $G$.

\textbf{Proof.} The corresponding Lie bracket on  ${\cal G}=I\oplus \mathbb{R}e$, where $I$ is an Abelian ideal, is given by $ [x,e]=u(x)$ for $x\in I$ and $u\in End(I)$. Consider the linear map $L:{\cal G}\rightarrow End(\cal G)$ defined by $L_{x}=0$ for $x\in I$ and $L_{e}=ad_{e}$. It is straightforward to check that $L_{a}b-L_{b}a=[a,b]$ and $L_{[a,b]}=[L_{a},L_{b}]$ for all $a$ and $b$ in ${\cal G}$. So $ab:=L_{a}b$ is a left symmetric product compatible with the bracket of ${\cal G}$. Moreover for any scalar $2$-cocycle $\omega$ of the Lie algebra ${\cal G}$ we have the formula 
\begin{center}
$\omega(L_{a}b,c)+ \omega( b,\omega L_{a}c)=0$.
\end{center}
Hence for $a^{+}$ and $b^{+}$, left invariant vector fields in $G$, the formula $\nabla_{a^+}b^{+}=(ab)^{+}$ , where $\omega$ is nondegenerate, determines a left invariant connection verifying the required conditions. \hfill $\square$\\

Some symplectic Lie groups or their dual or their doubles corresponding to some solutions of classical Yang-Baxter equations are endowed with other left invariant geometric structures for example flat semi-Riemannian metrics, complex structures, polynomial Poisson structures etc (see [12],[21],[22],[23]).

If the (simply connected ) symplectic Lie group $(G,\omega^{+})$ is quadratic i.e. has a bi-invariant metric given by a non degenerate bilinear form $k$ on ${\cal G}$, the affine structure defined by (1) is the Levi-Civita connection corresponding to the left invariant semi-Riemannian metric given by the bilinear form $<x,y>:=k(Dx,Dy)$ with $\omega(x,y)=k(Dx,y)$. In this case, the Lie algebra ${\cal G}$ is nilpotent because $D$ is an invertible derivation. Consequently any  nilmanifold $\Gamma\backslash G$ is endowed with a flat and complete semi-Riemannian metric. For more details ( see [7],[21]).

Notice that the inverse $r:=\omega^{-1}$ is an invertible solution of the classical Yang-Baxter  equation on ${\cal G}$. In [12]  it is shown that the Poisson tensors $r^{+}$ and  $r^{-}$ are polynomials. In fact  the Poisson-Lie tensor $\Pi(r)=r^{+}-r^{-}$ is polynomial of degree at most 2 (see [22]). Moreover, any  double Lie group $D(G)$ of $(G,\Pi(r))$ is a Manin Lie group endowed with an affine left invariant structure $\nabla$ and a left invariant complex structure $J$ such that $\nabla J=0$  (see [22] ). In these terms and as consequence of these results we have:

\textbf{Proposition 4. } Let $(G,\omega^{+})$ be a real simply connected nilpotent symplectic Lie group and $D(G)$ the double simply connected Lie group of $(G,\Pi(r))$. Then any lattice $\Gamma$ in $G$ determines naturally a class of lattices in $D(G)$ and if $\Gamma'$ is a lattice in $D(G)$ the  nilmanifold $\Gamma'\backslash D(G)$ is endowed with an homogeneous semi-Riemannian structure, two transversal  Lagrangians foliations, a flat complete connection $\nabla$ and a complex structure $J$ such that $\nabla J=0$.

\textbf{Proof.} Let $T^{*}G$ the cotangent bundle of $G$ endowed with the Lie structure semidirect produit of $G$ by ${\cal G}$ via the coadjoint action and recall that the map $\theta: D( {\cal G},\Pi(r))\rightarrow t^{*}\cal G $ with $\theta(\alpha,x):=(\alpha,r(\alpha)+x)$ is an isomorphism of Manin algebras ( see for example [22]) for any solution $r$ of the classical Yang-Baxter equation on $r$.

We will prove that $\Gamma$ determines lattices in $D(G)$. Now, the $\mathbb{Z}$-span ${\cal L}$ of the set $ exp_{G}^{-1}(\Gamma)$ is a lattice in the additive group ${\cal G}$. Consider a basis $\mathbb{B}$ of ${\cal G}$ contained in ${\cal L}$. It is clear that the set $\mathbb{B}':=(\mathbb{B}^{*}\times{0})\times ({0}\times\mathbb{B})$ where $\mathbb{B}^{*}$ is the dual basis of $\mathbb{B}$ is a basis of the Lie algebra $t^{*}\cal G={\cal G}^{*}\times{\cal G}$ with rational constants structure. Take now  the $\mathbb{Q}$-vector space generated by $\mathbb{B}'$, then if $L$ is any lattice of maximal rank of $t^{*}\cal G$ contained in this $\mathbb{Q}$ space, the group generated by $exp_{t^{*}\cal G}(L)$ is a lattice in $D(G)$. \hfill $\square$\\

\textbf{Acknowledgements.} We are grateful to the referee for helpful comments.













 



\end{document}